\documentclass[dvips,aos,preprint]{imsart}

\RequirePackage[OT1]{fontenc}

\usepackage{amsmath,amssymb,amsthm, latexsym,amsfonts}
\usepackage[dvips]{graphics,epsfig}

\newcommand{\R}{\mathbb{R}}
\newcommand{\N}{\mathbb{N}}
\newcommand{\B}{\mathbb{B}}

\newcommand{\FF}{\mathcal{F}}
\newcommand{\BB}{\mathcal{B}}
\newcommand{\CC}{\mathcal{C}}
\newcommand{\DD}{\mathcal{D}}

\newtheorem{theorem}{Theorem}[section]
\newtheorem{corollary}[theorem]{Corollary}
\newtheorem{remark}[theorem]{Remark}
\newtheorem{example}[theorem]{Example}

\newtheorem{lemma}[theorem]{Lemma}
\newtheorem{definition}[theorem]{Definition}
\numberwithin{equation}{section}

\newenvironment{proofof}{\noindent\sc Proof of}{
    \hspace*{\fill} $\Box$ \vspace{2ex} }

\newcommand{\Ind}[1]{{\boldsymbol{1}_{\textstyle\{#1\}}}}
\renewcommand{\theenumi}{\roman{enumi}}

\def\eps{\varepsilon}

\def\Min(#1,#2){#1\wedge #2}
\def\Max(#1,#2){#1\vee #2}
\newcommand{\ceil}[1]{\lceil #1 \rceil}
\newcommand{\floor}[1]{\lfloor #1 \rfloor}

\def\Cov{{\rm Cov}}

\def\e{{\rm e}}

\def\rueck{\noindent\hangafter=1 \hangindent=1.3em}

\arxiv{0910.0343}

\allowdisplaybreaks

\begin{document}

\begin{frontmatter}

\title{Limit Theorems for Empirical Processes of Cluster Functionals\protect\thanksref{T1}}
\runtitle{Empirical Cluster Processes}
\thankstext{T1}{We would like to thank Johan Segers
and Jon Wellner for very helpful discussions. We also want to thank
two referees and the editors for useful comments and in particular
for the suggestion to add a bootstrap example. The research work was
partly supported by
 the Swedish Foundation for Strategic Research through the Gothenburg Mathematical Modeling
 Centre.}

\begin{aug}
  \author{\fnms{Holger}  \snm{Drees}\corref{}\ead[label=e1]{holger.drees@uni-hamburg.de}},
  \and \author{\fnms{Holger} \snm{Rootz\'{e}n}\ead[label=e2]{rootzen@math.chalmers.se}}

  \affiliation{University of Hamburg and Chalmers and Gothenburg University}

  \address{University of Hamburg\\ Department of Mathematics,
SPST\\ Bundesstr.\ 55\\ 20146 Hamburg\\ Germany\\
          \printead{e1}}

  \address{Chalmers University \\ Department of Mathematical Statistics\\
412 96 G\"{o}teborg\\ Sweden\\ and\\
 Gothenburg University\\ Department
of
Mathematical Statistics\\  412 96 G\"{o}teborg\\ Sweden\\
          \printead{e2}}

\end{aug}

\begin{abstract}
Let $(X_{n,i})_{1\le i\le n, n\in\N}$ be a triangular array of
row-wise stationary  $\R^d$-valued random variables. We use a
``blocks method'' to define clusters  of extreme values: the rows of
$(X_{n, i})$ are divided into $m_n$ blocks $(Y_{n,j})$, and if a
block contains at least one extreme value the block is considered to
contain a cluster. The cluster  starts at the first extreme value in
the block and ends at the last one. The main results are uniform
central limit theorems for empirical processes $ Z_n(f) := \frac
1{\sqrt{n v_n}} \sum_{j=1}^{m_n} \big( f(Y_{n,j})- E
f(Y_{n,j})\big), $ for $v_n=P\{X_{n, i} \neq 0\}$ and $f$ belonging
to classes of cluster functionals, i.e.\ functions of the blocks
$Y_{n,j}$ which only depend on the cluster values and which are
equal to 0 if $Y_{n,j}$ does not contain a cluster. Conditions for
finite-dimensional convergence include $\beta$-mixing, suitable
Lindeberg conditions and convergence of covariances. To obtain full
uniform convergence we use either ``bracketing entropy''  or bounds
on covering numbers with respect to a random semi-metric.  The
latter makes it possible to bring the powerful
Vapnik-\v{C}ervonenkis theory to bear. Applications include
multivariate tail empirical processes and empirical processes of
cluster values and of order statistics in clusters. Although our
main field of applications is the analysis of extreme values, the
 theory can be applied more generally to rare events
occurring e.g.\ in nonparametric curve estimation.
\end{abstract}



\end{frontmatter}

\section{Introduction}

The next challenge for extreme value statistics is modeling and
estimation of the structure of clusters of extreme values. As one
concrete example, the Europe 2003 heat wave may have killed around
60,000 persons. There has been a substantial discussion of whether
it could be attributed to global warming. The Nature paper Stott
{\em et al.} (2004) uses extreme value methods with average summer
temperature as a proxy for a heat wave to try to answer this
question. However, the health effects are in reality linked to clusters of
extremely high temperatures over much shorter time periods, and the
fluctuations of temperature during this period determine risks.

Similarly, river flooding may be caused by not just one extreme
rainfall event, but also by the ground already being saturated with
water due to high precipitation during the preceding 5-10 days. This
was e.g. the case for the large flood which occurred in Northern
Sweden on July 26, 2000. Thus, again, an entire sequence of large
values are at the center of interest.

This paper develops an empirical limit theory for clusters of
extremes in stationary sequences. It provides a unified basis for
asymptotic analysis of statistical methods which aim at answering
questions such as the ones above. Results include limit theorems for
tail array sums, in particular for multivariate tail empirical
processes, and for joint survival functions of the values and order
statistics in a cluster. More special examples such as upcrossings,
compound insurance claims,  kernel density and bootstrap
estimators, are also studied.

Estimation of the extremal index (roughly, the inverse of the
expected clusters length) has received substantial attention in the
extreme value statistics literature. The results of this paper can
be used to prove asymptotic normality for a general type of
estimators based on blocks of exceedances, see Drees (2009). There
are also a few papers (e.g.\ Bortot and Tawn (1998), Sisson and
Coles (2003)) on Markov chain modeling of clusters of extreme
values. However, a major part of the work to develop useful
statistical methods for the structure of clusters of extremes still
remains to be done. Our goal is that this paper will be useful for
the analysis of existing methods, and that it will spur development
of new methods.

More specifically, we consider triangular arrays of row-wise
stationary sequences of random variables. The variables are assumed
to take their values in some set $E \subset \R^d$, with $E=\R$ and
$E=\R^d$ as the standard examples. Clusters of extremes are defined
through a ``blocks'' method. The variables in each row of the array
are divided up into blocks, and a cluster of extremes starts with
the first ``extreme'' value in a block, if there is such a value,
and ends with the last one. Such a cluster is termed the ``core'' of
the block. A function which maps a block into a real number is
called a ``cluster functional'' if it only depends on the core of
the block and if it equals 0 for blocks without extremes. In
contrast to standard uniform central limit theorems, cores (i.e.\
clusters of extremes) consist of a random number of variables, and
hence cluster functionals have to be defined on a space of vectors
of arbitrary lengths.

The aim is to prove uniform central limit theorems for interesting
classes of cluster functionals. We throughout use $\beta$-mixing
(or, with another name, absolute regularity) as the basic dependence
restriction. It is very widely applicable and makes it possible to transfer calculations from
dependent blocks to easier calculations with independent blocks.
Finite-dimensional convergence of the cluster functionals in
addition requires Lindeberg conditions and convergence of
covariances. We use suitable formulations of ``bracketing entropy''
to give conditions for asymptotic tightness, and bounds on covering
numbers with respect to a random semi-metric to prove asymptotic
equicontinuity. The latter in particular makes it possible to use
Vapnik-\v{C}ervonenkis theory to prove asymptotic equicontinuity. As
usual uniform central limit theorems follow from finite-dimensional
convergence together with asymptotic tightness, or together with
asymptotic equicontinuity.

In the important context of estimation for panel count data, two
articles by Wellner and Zhang (2000, 2008) use uniform central limit
theory for vectors of random lengths. These articles are aimed at
the specific application and not at general theory. Hence they use
special properties (such as monotonicity) of the classes of
functions, do not consider triangular arrays, assume that the
vectors are independent, and, in the second paper also that the
lengths of the vectors are uniformly bounded. However,  the basic
tools to prove tightness, i.e.\ random covering numbers for the
general case, and bracketing entropy for the uniformly bounded case
are the same as in the present paper.
 We have not found any other references
on uniform central limit theory for random vectors with random
lengths.

One application of the theory of this paper is to multivariate tail
empirical processes for stationary time series. Let $(X_i)_{i\in\N}$
be a time series with marginal survival function $\bar H=1-H$. The
univariate tail empirical process is defined as
$$ e_n(x) := \frac 1{\sqrt{nv_n}} \sum_{i=1}^n \big( \Ind{X_{n,i}>x}-\bar
H(u_n+a_nx)\big), \quad x\in[0,\infty),
$$
where
\begin{equation}  \label{standexceed}
 X_{n,i}:= \Big( \frac{X_i-u_n}{a_n}\Big)_+ = \max \Big(
\frac{X_i-u_n}{a_n},
    0\Big), \quad 1\le i\le n.
\end{equation}
The multivariate tail empirical process is defined analogously, see
Examples \ref{tailedfexam} and \ref{tailedfexam2} below. In the
definition $(u_n)_{n\in\N}$ is an increasing sequence of thresholds
such that $v_n := P\{X_1>u_n\}\to 0$, and $(a_n)_{n\in\N}$ is a
sequence of positive normalizing constants such that the conditional
distribution of $X_{n,1}$ given that $X_{n,1}>0$ converges weakly to
some non-degenerate limit. (In particular, the distribution function
(df) of $X_1$  then belongs to the domain of attraction of some
extreme value distribution.) Rootz\'{e}n (1995, 2009) proved weak
convergence of $e_n$ to a Gaussian process; see Example
\ref{tailedfexam2} for details. Such limit theorems have proved
quite useful for semi-parametric statistical analysis of the
marginal tail behavior (Drees, 2000, 2002, 2003). The present paper
extends convergence to multivariate tail empirical processes and
makes a small improvement of the results in Rootz\'en (2009).

Tail empirical processes do not capture information on location in
the extreme clusters, and hence do not catch the serial extremal
dependence structures which are at the center of interest in
connection with e.g.\ heat waves or river floods. A second class of
applications of our main theorems is to joint survival functions and
joint distributions of the order statistic of the values within an
extreme cluster.

The paper is organized as follows. In Section 2 we first introduce
empirical processes of cluster functionals. This generalizes
concepts first introduced by Yun (2000) and developed further by
Segers (2003). We then derive uniform central limit theorems for
these empirical processes under quite general abstract conditions.
Sections 3 contains applications to tail array sums, with the
multivariate tail empirical process as a prominent example. In
Section 4 we consider empirical processes of indicator variables,
and in particular joint distributions of variables and of the order
statistics in the clusters of extreme values. Proofs are given in
Section 5.

\section{Limit theorems for general empirical cluster processes}
\label{sec:limittheorems}

This section first sets out the basic definitions and assumptions
which are used throughout the paper and then, in Subsection
\ref{sect:fidi}, gives conditions for finite-dimensional convergence
of the empirical processes $(Z_n(f))_{f\in\FF}$ (defined below). The
following subsections consider asymptotic tightness and asymptotic
equicontinuity of these empirical processes. As usual,
finite-dimensional convergence together with either asymptotic
tightness or asymptotic equicontinuity  gives convergence of $Z_n$
in the space $\ell^\infty(\FF)$ of bounded functions indexed by
$\FF$.

For some $d\in\N$, let $E$ be a measurable subset of $\R^d$
containing  0 and let $(X_{n,i})_{1\le i\le n, n\in\N}$ be a
triangular array of row-wise stationary random variables (rv's) with
values in $E$. Typically the $(X_{n, i})$ have been obtained by
``renormalization'' of some other process, where the renormalization
maps all non-extreme values to 0. A generic example (cf.\ the
introduction) is $E=\R$ and  $X_{n, i} = (\frac{X_i - u_n}{a_n})_+$
where $(X_i)_{i\in\N}$ is a stationary univariate time series. Here
$u_n $ tends to the right endpoint of the support of $X_i$, so that
$X_{n, i}$ is 0 unless $X_i$ is ``large'', i.e.\ unless $X_i > u_n$.

The ``empirical process $Z_n$ of cluster functionals'' is defined as
$$ Z_n(f) := \frac 1{\sqrt{n v_n}} \sum_{j=1}^{m_n} \big( f(Y_{n,j})-
E f(Y_{n,j})\big), \quad f\in\FF.
$$
Here $Y_{n,j}$ is the $j$-th block of $r_n$ consecutive values of
the $n$-th row of $(X_{n,i})$. Thus there are
$m_n:=\floor{n/r_n}:=\max\{j\in\N_0\mid j\le n/r_n\}$ blocks
$$ Y_{n,j} := (X_{n,i})_{(j-1)r_n+1\le i\le jr_n}, \quad 1\le j\le
m_n,
$$
of length $r_n$. We write $Y_n$ for a ``generic block'' so that $Y_n
\stackrel{d}{=} Y_{n,1}$. The block lengths $r_n$ tend to infinity,
but slower than $n$, and
$$ v_n := P\{X_{n,1}\ne 0\} \to 0.$$ Further $\FF$ is a class of ``cluster functionals'',
i.e.\ functions which only depend on the part of the block which
contains all nonvanishing observations, see below.

In the univariate case $E=\R$, cluster functionals have been
introduced by Yun (2000) and Segers (2003). The definition is as
follows.
\begin{definition}
  \begin{enumerate}
    \item The set $E_\cup := \bigcup_{l\in\N} E^l$ of vectors of arbitrary length is
    equipped with the $\sigma$-field $\mathbb{E}_\cup$ that is induced by
    the Borel-$\sigma$-fields on $E^l$, $l\in\N$.
    \item For an arbitrary $k\in\N$ and $x=(x_1,\ldots,
    x_k)\in E^k$ the {\em core}
    $x^c\in E_\cup$ of $x$ is defined by
    $$ x^c := \left\{
      \begin{array}{l@{\quad}l}
        (x_l)_{l_1\le l\le l_2} & \text{if } x\ne(0,\ldots,0),\\
        0 & \text{otherwise,}
      \end{array} \right.
    $$
    where
    \begin{eqnarray*}
      l_1 & := & \min\big\{ i\in\{1,\ldots,k\} \mid x_i\ne 0\big\}
      \\
      l_2 & := & \max\big\{ i\in\{1,\ldots,k\} \mid x_i\ne 0\big\}
    \end{eqnarray*}
    The {\em length of the core} of $x$ is defined as
    $L(x):=l_2-l_1+1$ if $x^c\ne 0$ and $L(x)=0$ if $x^c=0$.
    \item A measurable map $f:(E_\cup,\mathbb{E}_\cup)\to(\R,\B)$ is called
    a {\em cluster functional} if $f(x)=f(x^c)$ for all
    $x\in E_\cup$, and $f(0)=0$.
  \end{enumerate}
\end{definition}

Typical examples are functionals of the type
$$ f(x_1,\ldots,x_k) := \sum_{l=1}^k \phi(x_l)  $$
where $\phi:E\to \R$ satisfies $\phi(0)=0$, which are related to
so-called tail array sums, and, in the case $E=[0,\infty)$,
$$ f(x_1,\ldots, x_k):= \max_{1\le i\le k} x_i, $$
which corresponds to the (componentwise) maximum of a cluster. Many
more examples will be discussed in the Sections 3 and 4.

The proofs below will use the well-known ``big blocks, small
blocks''  technique together with a $\beta$-mixing condition to boil
down convergence to convergence of sums over i.i.d.\ blocks. The
$\beta$-mixing coefficients (also called the coefficients of
absolute regularity) for  $(X_{n,i})_{1\le i\le n}$  are defined by
$$ \beta_{n,k} := \sup_{1\le l\le n-k-1} E \Big(
\mathop{\text{sup}}_{B\in\BB_{n,l+k+1}^n} |
P(B|\BB_{n,1}^l)-P(B)|\Big) $$ where $\BB_{n,i}^j$ denotes the
$\sigma$-field generated by $(X_{n,l})_{i\le l\le j}$. Since the
$X_{n,i}$ take values in a Polish space, the supremum can be taken
over a countable set of $B$'s, and hence is measurable. (On general
spaces ``sup'' has to be replaced by ``ess-sup'', which is defined
as a measurable function which is a.s.\ larger than or equal to $|
P(B|\BB_{n,1}^l)-P(B)|$ for all $B\in\BB_{n,l+k+1}^n$ and a.s.\
smaller than or equal to all other measurable functions with this
property.) In addition to the $\beta$-mixing coefficients and the
lengths $r_n$ of the big blocks, the ``big blocks, small blocks''
technique uses an intermediate sequence $\ell_n$ of integers, the
lengths of small blocks which are used to separate the big blocks in
the proofs.

Throughout we will use the following Basic Assumptions.
\\[1ex]
{\bf (B1)}\; \; \; \parbox[t]{11.2cm}{The rows $(X_{n,i})_{1\le i\le
n}$ are stationary, $ \ell_n = o(r_n),   \; \ell_n\to\infty,$ \\
$r_n = o(n), \; r_nv_n \to 0, \;  nv_n \to \infty,$}
\\[1ex]
and
\\[1ex]
{\bf (B2)} \; \; \; $\beta_{n,l_n} \frac n{r_n} \to 0. $ \\[1ex]
Sometimes we will also use the assumption
\\[1ex]
{\bf (B3)} \; \; \;$
        \lim_{m\to\infty} \limsup_{n\to\infty} \beta_{n,m}  =  0.$
\vspace*{2mm}

It follows from $r_nv_n \to 0$ that $v_n \to 0$ and hence that non-zero values of $X_{n,i}$
are rare events. The most important example we have in mind are the
standardized excesses given in \eqref{standexceed}. However, other
examples occur in the context of nonparametric density estimation or
nonparametric regression in a natural way (cf.\ Example
\ref{densityexam}). Since $nv_n$ is the expected number of nonzero
values of $(X_{n,i})_{1 \leq i \leq n}$, the assumption $nv_n \to
\infty$ seems necessary if one wants to obtain normally distributed
limits.

More specifically, the assumption  $r_nv_n \to 0$ means that the
probability of a block being non-zero tends to zero. In particular,
it implies that if the  row variables are i.i.d., then
asymptotically cores -- or equivalently clusters of ``extremes'' --
will have length one, as they intuitively should have. To see this
note that if the variables in a row are independent, then
asymptotically the number of non-zero values in a block of length
$r_n$ has a Poisson distribution with mean $r_nv_n$ and that then
the conditional probability that there are more than one non-zero
value in a block, given that there is at least one non-zero value is
(approximately) $(1-e^{-r_nv_n} - r_nv_n e^{-r_nv_n})/(1-
e^{-r_nv_n})$. This tends to zero if and only if $r_nv_n \to 0$.

For a given sequence $(r_n)_{n\in\N}$,
Assumption (B2) requires a minimum rate at which the mixing
coefficients $\beta_{n,l}$ tend to 0 as $l\to\infty$. The condition
(B3) e.g.\ holds if the $X_{n,i}$ are obtained by renormalizing a
single absolutely regular process.

\begin{remark} \rm
\begin{enumerate}
\item
The proofs of Theorems \ref{fidisconvtheo} and \ref{brackettheo},
  of Lemma \ref{lemma:clusterconv} (ii) and (iii) and of Lemma \ref{iidcompprop} below
 in fact do not use the assumption $r_nv_n \to 0$ of (B1),
 but only that $v_n \to 0$. The same remark applies to Theorem \ref{condrem2}
 if one replaces condition (D5) below by the following slightly
 stronger version:
 For all $\delta>0, n\in\N, l\in\{0,1\}, (e_i)_{1\le i\le
   \floor{m_n/2}+1} \in \{-1,0,1\}^{\floor{m_n/2}+1}$ and $k\in\{1,2\}$ the
   map $\sup_{f,g\in\FF, \rho(f,g)<\delta}$
   $\sum_{j=1}^{\floor{m_n/2}+l} e_j\big(
   f(Y_{n,j}^*)-g(Y_{n,j}^*)\big)^k$
      is measurable.

 Hence these results hold also if
 the assumption $r_nv_n \to 0$ is replaced by the  weaker $v_n \to 0$.
\item
  It is not essential that $E$ is a subset of $\R^d$. Indeed, one may assume
  that $X_{n,i}$ takes on values in an arbitrary set $E$. Then
  one chooses some special element $e_0\in E$ which takes over the
  role of 0. In this more general setting, a cluster functional
  is defined as a functional on $\bigcup_{l\in\N} E^l$ whose value is not
  changed if $e_0$ is added at the beginning or at the end of some
  vector in $\bigcup_{l\in\N} E^l$.
  \hspace*{\fill} $\Box$
  \end{enumerate}
\end{remark}

\subsection{Convergence of fidis}\label{sect:fidi}

We first give a general result on the convergence of the
finite-dimensional marginal distributions (fidis), and then
introduce simpler, but more restrictive assumptions, which also are
sufficient for convergence. Proofs are deferred to Section
\ref{proofs}.

We will use the notation $x^{(k)}$ for the vector $(x_1,
\ldots,x_k)$  made up by the first $k$ components in the vector $x$,
if $x$ has at least $k$ components, and otherwise $x^{(k)}=x$.
Similarly we write $x^{(\ell; k)} = (x_\ell, \ldots, x_k)$ for the
vector consisting of components number $\ell$ to number $k$ in $x$,
if $x$ has at least $k$ components, and otherwise $x^{(\ell; k)}$
starts at component no.~$\ell$ and ends at the end of $x$ (if $x$ is
shorter than $\ell$ then $x^{(\ell; k)} = 0$). As before let $\FF$
be a class of cluster functionals, recall that $Y_n \stackrel{d}{=}
Y_{n,1}$, where $Y_{n,1}$ is the first block in the $n$-th row. For $f \in \FF$ write
$$
 \Delta_n(f)  :=   f(Y_n) - f(Y_n^{(r_n - \ell_n)})
$$
for the difference between $f$ evaluated at the $r_n$ components of
the entire block and $f$ evaluated at the first $r_n - \ell_n$
components of the  block. The general ``Convergence Conditions'' are
as follows.
\\[1ex]
{\bf (C1)}\vspace*{-0.6cm}

\begin{eqnarray}
      E\Big( (\Delta_n(f)- E\Delta_n(f))^2
      \Ind{|\Delta_n(f)- E\Delta_n(f)|\le
      \sqrt{nv_n}}\Big) & = & o(r_n v_n)   \label{necsuffbbsbcond1} \nonumber\\
       E\Big( (\Delta_n(f)- E\Delta_n(f))
      \Ind{|\Delta_n(f)- E\Delta_n(f)|>
      \sqrt{nv_n}}\Big) & = & o\Big(r_n\sqrt{\frac{v_n}{n}}\Big)   \nonumber\\ 
     P \big\{|\Delta_n(f)- E\Delta_n(f)|>
     \sqrt{nv_n}\big\} & = & o(r_n/n) \label{necsuffbbsbcond2}
     \nonumber
   \end{eqnarray}
   for all $f\in\FF$.
\\[1ex] {\bf(C2)}\; \; \; \parbox[t]{11.2cm}{
$E\Big( (f(Y_n) - E f(Y_n))^2\Ind{|f(Y_n) - E
    f(Y_n)|>\eps\sqrt{nv_n}}\Big) = o(r_n v_n),$ \\\hspace*{8cm} $\forall\,
    \eps>0, f\in\FF.$}
\vspace*{2mm}
\\[1ex] {\bf (C3)}\; \; \;  $\displaystyle \frac 1{r_nv_n} Cov\big( f(Y_n),
g(Y_n)\big) \to c(f,g) \quad \forall\, f,g\in\FF.$

\vspace*{2mm}

The block $Y_n^{(r_n - \ell_n)}$ is obtained from $Y_n$ by omitting
a small block of $l_n$ observations at the end. Accordingly (C1)
means that asymptotically this omission does not influence the fidis
of the empirical process of cluster functionals (see the proof of Lemma
\ref{iidcompprop}). By the definition
of cluster functionals this is usually fulfilled if with high
probability there are few or no non-zero observations in the omitted
short blocks. Specifically, if components number $r_n-l_n+1\le i\le
r_n$ all are zero, then $Y_n$ and $Y_n^{(r_n - \ell_n)}$ have the
same core, and thus $\Delta_n(f)=0$.

Assumption (C2) is the standard Lindeberg condition. The assumption
of convergence of covariances, (C3), is the final ingredient needed
to ensure finite-dimensional convergence in the present triangular
array setup.

\begin{theorem} \label{fidisconvtheo}
  Suppose the basic assumptions (B1) and (B2) hold, and that (C1)--(C3) are satisfied.
   Then the fidis of the empirical process $(Z_n(f))_{f\in\FF}$ of cluster
  functionals converge to the fidis of a Gaussian process $(Z(f))_{f\in\FF}$
  with covariance function $c$.
\end{theorem}

In general, the convergence (C3) of the covariance function must be
verified directly. However, we also give additional sufficient
conditions which are simpler to verify in some situations. A first
very simple version, (C3'), requires convergence only after
``truncation'' to a fixed (but arbitrary) length. Before stating it
we recall the notation $L(Y_n)$ for the length of the core of $Y_n$.
\\[1ex]{\bf (C3')}\; \; \; \parbox[t]{10.8cm}{
For  $f \in \FF$ it holds that
\begin{equation}  \label{eq:C3.1'}
\lim_{k\to\infty}\limsup_{n\to\infty} \frac{1}{r_nv_n}
E\big(f(Y_n)^2 \boldsymbol{1}_{\{L(Y_n)
> k\}}\big) = 0,
\end{equation}
 and for $f, g \in \FF$ there is a sequence $R_{n,k}$ with $\lim_{k \to
\infty} \limsup_{n \to \infty} |R_{n,k}| = 0$ such that
\begin{equation}  \label{eq:C3.2'}
  \lim_{n\to\infty}  \frac{1}{r_nv_n}
   \big(E\big(f(Y_n)g(Y_n) \boldsymbol{1}_{\{L(Y_n) \leq k\}}\big) + R_{n,k}  = c_k(f,g).
\end{equation}}

A typical situation when \eqref{eq:C3.1'} holds is when the cluster
lengths $(L(Y_n))_{n=1}^\infty$ are tight under $P(\cdot|\; Y_n \neq
0)$ and $\big(f(Y_n)^2\big)_{n\in\N}$ is uniformly integrable under
$P(\cdot\mid Y_n \neq 0),$ for $f \in \FF$. This follows from the
observation that $\frac{1}{r_nv_n}|E( \cdot)| \leq |E\big( \cdot
\mid Y_n \neq 0 \big)|$,  which in turn follows from $P(Y_n \neq 0)
\leq r_nv_n$.

 In a second  assumption (C3'') we generalize the powerful results of Segers (2003) to the present
abstract setting.  In doing this we do not aim at the greatest
possible generality, but give versions which suit our purposes best.
It may be noted that unlike in the situation considered by Segers,
in general weak convergence of the indicators
$\boldsymbol{1}_{\{0\}}(X_{n,i})$ does not follow from weak
convergence of $X_{n,i}$.  In the statement of the condition we use
that the value of a cluster functional $f$ applied to a sequence
$(x_i)_{i\in\N}$ with $m_x:=\sup\{i\in\N\mid x_i\ne 0\}<\infty$ can
be defined in a natural way as $f((x_i)_{1\le i\le m_x})$. The
conditions are as follows.
\bigskip

{\bf (C3'')}\parbox[t]{11cm}{
\begin{itemize}
  \item[(C3.1'')] There is a sequence $W=(W_i)_{i\in\N}$ of $E$-valued r.v.'s
such that for all $k \in \N$, the joint conditional distribution
$P^{(X_{n,i}, \boldsymbol{1}_{\{0\}}(X_{n,i}))_{1\le i\le k}\mid
X_{n,1}\ne 0}$ converges weakly to
$P^{(W_i,\boldsymbol{1}_{\{0\}}(W_i))_{1\le i\le k}}$, and all $f\in
\FF$ are a.s.\ continuous with respect to the distributions of
$W^{(k)}$ and $W^{(2;k)}$, for  all $k$, i.e.
    \begin{eqnarray}  \label{fcontcond}
      \lefteqn{P\{W^{(2; k)}\in D_{f,k-1}, W_i=0\;\forall\, i>k\}} \nonumber \\
      & = &
      P\{W^{(k)}\in D_{f,k}, W_i=0\;\forall\, i>k\}=0,
    \end{eqnarray}
    with $D_{f,k}$ denoting the set
    of discontinuity points of $f_{|E^k}$.

  \item[(C3.2'')] For all $f \in \FF$ the sequence  $\big(f(Y_n)^2\big)_{n\in\N}$
  is uniformly integrable under $P(\cdot)/(r_nv_n)$.
\end{itemize}}
 \vspace*{2mm}

Again, (C3.2'') is implied by the perhaps more intuitive condition that
$\big(f(Y_n)^2\big)_{n\in\N}$ is uniformly integrable under
$P(\cdot\mid Y_n \neq 0)$.

In the proof of the next two results we in fact will use a slightly weaker
(but instead more complicated) version of \eqref{fcontcond}, see Remark \ref{condC3ppgen} below.

\begin{corollary}  \label{covcondcorol}
  Suppose that (B1), (B2), and (C1) are satisfied. If furthermore either (C2) and (C3')
  or else (B3) and (C3'')
  hold,  then
  the fidis of the empirical process $(Z_n(f))_{f\in\FF}$ of cluster
  functionals converge to the fidis of a Gaussian process
  $(Z(f))_{f\in\FF}$. Specifically, (C3') implies that (C3) holds and
  that the covariance function $c$ of $Z$ is
  obtained as
  $$
  c(f,g) = \lim_{k \to \infty} c_k(f,g).
  $$
If (C3'') holds, then
\begin{equation} \label{thetacovconv}
c(f,g) = E\big((fg)(W) - (fg)(W^{(2; \infty)})\big).
\end{equation}
\hfill $\Box$
\end{corollary}

Equation \eqref{thetacovconv} is explained in Lemma
\ref{lemma:clusterconv} below. It generalizes the most important
results of Segers (2003) to the present more abstract setting.

\begin{lemma}  \label{lemma:clusterconv}
  \begin{enumerate}
    \item If (B1) and (B3) hold,
      then
      \begin{equation}  \label{clusterfctnalapprox}
        E\big( f(Y_n) \mid Y_n\ne 0\big) = \frac 1{\theta_n}
        E\Big( f\big( X_{n}^{(r_n)}\big) - f\big( X_{n}^{(2, r_n)}\big) \mid X_{n,1}\ne 0\Big) + o(1)
      \end{equation}
      where the term $o(1)$ tends to 0 as $n$ tends to $\infty$ uniformly for all cluster
      functionals $f$ such that $\|f\|_\infty\le C$,  for any $C\in\R$, and
     $$
       \theta_n  :=  \frac{P\{Y_n\ne 0\}}{r_n v_n}
          =   P\big(      X_{n}^{(2;r_n)}=0 \mid X_{n,1}\ne 0\big)      (1+o(1)).
      $$
    \item If (B1), (B3), and the  assumption of (C3.1'')
    all are satisfied, then
    \begin{equation} \label{eq:finitew}
    m_W=\sup\{i \geq 1 \; | \; W_i \neq 0\} < \infty
    \end{equation}
    and
      $$
       \lim_{n\to\infty} \theta_n
       =  \theta := P\{W_i=0\; \forall\, i\ge 2\}=P\{m_W=1\}>0.
     $$

    \item   If (B1), (B3), and (C3.1'') hold, then the conditional distribution\\
    $P^{f(Y_n)\mid Y_n\ne 0}$ converges weakly to the
    probability measure
    $$ \mu_{f,W} := \frac 1{\theta} \Big(
    P\big\{f(W)\in \cdot \big\}- P\big\{f\big(W^{(2;\infty)}\big)\in \cdot, m_W\ge 2 \big\}
    \Big).
    $$
  \end{enumerate}
\end{lemma}

Note that $\mu_{f,W}(\R)=1$ by (ii). However, it is not so obvious
that $\mu_{f,W}$ is indeed a positive (and hence a probability)
measure.

\begin{remark}  \label{condC3ppgen} \rm
  We will prove Corollary \ref{covcondcorol}
  and Lemma \ref{lemma:clusterconv} under the following weaker version of the continuity
  assumption \eqref{fcontcond}:\\[0.5ex]
  For $k\in\N$ and $I\subset\{1,\ldots,k\}$ let $N_{k,I}:=\{x\in
E^k\mid x_i=0,\,\forall i\in I,\; x_i\ne 0,\, \forall i\not\in I\}$
and denote by $D_{f,k,I}$ the set of discontinuity points of
$f|_{N_{k,I}}$. Then we assume
    \begin{eqnarray}
        & & P\{W^{(k)}\in D_{f,k,I}, W^{(k+1,\infty)}=0\}=0,\;
         \forall\; k\in\N, I\subset\{1,\ldots,k\},   \label{fcontcond1}\\
        & & P\{W^{(2;k)}\in D_{f,k-1,I}, W^{(k+1,\infty)}=0\}=0,\; \forall\; k\ge 2, I\subset\{1,\ldots,k-1\}. \label{fcontcond2}
    \end{eqnarray}
  This version can be used in some examples where \eqref{fcontcond}
  is not satisfied, because the boundary of $[0,\infty)^k$ belongs to the discontinuity sets
  $D_{f,k}$ and, according to Lemma \ref{lemma:clusterconv} (ii),
   the rv $W_i$ equals 0 with positive probability for
  $i>1$. \hfill$\Box$
\end{remark}

In the situation considered by Segers (2003) (i.e.\ with $X_{n,i}$
defined by \eqref{standexceed} for a stationary time series whose
finite-dimensional marginal distributions all belong to the domain
of attraction of some extreme value distribution), the sequence
$(W_i)_{i\in\N}$ is related to the so-called tail sequence (or tail
chain) $(U_i)_{i\in\N}$ (cf.\ Segers, 2003, Theorem 2) via
$W_i=\max(U_i,0)$. Then (C3'') is automatically satisfied, e.g., for
bounded cluster functionals if $D_{f,m}$ is a Lebesgue null subset
of $(0,\infty)^m$ for all $m$ and $f\in\FF$, because the rv's $U_i$
are continuous.

Further simpler, but more restrictive, sufficient conditions are
given in Lem\-ma~\ref{simplerfidicond} below. In particular, for
bounded cluster functionals one obtains
\begin{corollary} \label{bddfunctfidicond}
  If $\|f\|_\infty=\sup_{x\in E_\cup} |f(x)|<\infty$ for all
  $f\in\FF$ and
  the conditions (B1), (B2), (B3) and (C3.1'') hold,
  then the fidis of the empirical process $(Z_n(f))_{f\in\FF}$ of cluster
  functionals converge to the fidis of a Gaussian process
  $(Z(f))_{f\in\FF}$ with covariance function $c$ defined by
  \eqref{thetacovconv}.
\end{corollary}

\subsection{Asymptotic tightness}

In this subsection we give conditions which ensure asymptotic
tightness of $Z_n$ in the space $\ell^\infty(\FF)$. As a consequence
uniform central limit theorems for $Z_n$ hold if in addition the conditions of
Theorem \ref{fidisconvtheo} are satisfied. The alternative route via
asymptotic equicontinuity is considered in the next subsection.

In general the supremum of $Z_n(f)$ taken over uncountably many
cluster functionals $f$ need not be measurable. Hence, in some
instances, one has to work with outer probabilities and
expectations, denoted by $P^*$ and $E^*$ in the following; see van
der Vaart and Wellner (1996), Section 1.2, for details. The sequence
$(Z_n)_{n\in\N}$ is asymptotically tight if to any $\epsilon >0$
there is a compact set $K \subset  \ell^\infty(\FF)$ such that
$$
\limsup_{n \to \infty} P^*(Z_n \notin K^\delta) < \epsilon,   \;\;\; \mbox{for any} \;\;\; \delta > 0.
$$
Here $K^\delta$  is the set of elements in $\ell^\infty(\FF)$ which
are at most  a distance $\delta$ away from $K$.

We will use the assumptions (D1)--(D4) below  to prove tightness.
The first two assumptions  in various ways restrict the sizes of the
functions in $\FF$. In particular (D1) ensures that sample paths of
$Z_n$ belong to the space  $\ell^\infty(\FF)$ of bounded functions
on $\FF$. The assumption (D3) is an asymptotic continuity condition
on the covariance function which  is needed to ensure that the
limiting process has continuous sample paths. The most crucial
condition, (D4), restricts the complexity of the index set $\FF$ via
the so-called bracketing entropy. To state this assumption, the
following concept is needed.

The {\em bracketing number} $N_{[\cdot]}(\eps,\FF,L_2^n)$ here is defined
as the smallest number $N_\eps$ such that for each $n\in\N$ there
exists a partition $(\FF_{n,k}^\eps)_{1\le k\le N_\eps}$ of $\FF$
such that
\begin{equation}
  E^* \sup_{f,g\in\FF_{n,k}^\eps} \big(f(Y_{n})-g(Y_{n})\big)^2
  \le \eps^2 r_n v_n, \quad \forall\, 1\le k\le N_\eps.
\end{equation}

The assumptions are as follows.
 \vspace*{2mm}
\\[1ex]{\bf (D1)}\;\;\; \parbox[t]{11.3cm}{The index set $\FF$ consists of cluster
functionals $f$ such that $E(f(Y_n)^2)$ is finite for all $n\geq 1$
and such that the envelope function
$$ F(x) := \sup_{f\in\FF} |f(x)| $$
is finite for all $x\in E_\cup$.} \vspace*{2mm}
\\[1ex]{\bf (D2)}  $$E^* \Big( F(Y_{n}) \Ind{F(Y_{n})>\eps\sqrt{nv_n}}\Big)  =
    o\big(r_n\sqrt{v_n/n}\big), \quad \forall\, \eps>0.$$
\vspace*{2mm}
\\[1ex]{\bf (D3)}\;\;\;  \parbox[t]{11.1cm}{There exists a semi-metric $\rho$ on
$\FF$ such that $\FF$ is totally bounded (i.e., for all $\eps>0$ the
set $\FF$ can be covered by finitely many balls with radius $\eps$
w.r.t.\ $\rho$) such that  $$\lim_{\delta\downarrow 0}
\limsup_{n\to\infty}
    \sup_{f,g\in\FF, \; \rho(f,g)<\delta} \frac 1{r_nv_n}
    E(f(Y_{n})-g(Y_{n}))^2  = 0.$$}\\[1ex]
\vspace*{2mm}
\\[1ex]{\bf (D4)} \parbox[t]{11.1cm}{$$\lim_{\delta\downarrow 0} \limsup_{n\to\infty} \int_0^\delta
     \sqrt{\log N_{[\cdot]}(\eps,\FF,L_2^n)}\, d\eps  =  0.$$}

\begin{theorem}  \label{brackettheo}
  If the basic assumptions (B1) and (B2) hold and  (D1)--(D4) are
  satisfied, then the process $Z_n$ is asymptotically tight in $\ell^\infty(\FF)$.
If in addition the finite-dimensional distributions converge (which
in particular holds if (C1)--(C3) also are satisfied), then $Z_n$
converges to a Gaussian process $Z$ with covariance function
  $c$.
\end{theorem}

We collect a number of comments and variations of the conditions of
the theorem in the following remark. In particular we consider a
strengthened version (D2') of (D2).
 \vspace*{2mm}
\\[1ex]{\bf (D2')}\begin{equation*}
        E^* \Big( F^2(Y_{n}) \Ind{F(Y_{n})>\eps\sqrt{nv_n}}\Big)
        = o(r_nv_n), \quad\forall\, \eps>0.
    \end{equation*}

The proof of part (ii) of the remark is given in Section~\ref{proofs}

\begin{remark}  \label{condrem1}  \rm
  \begin{enumerate}
        \item If, for all $\eps>0$, there exists a partition $(\FF_k^\eps)_{1\le k\le
    N_\eps}$ of $\FF$ which does {\em not} depend on $n$ and which
    satisfies
    $$ E^* \sup_{f,g\in\FF_{k}^\eps} \big(f(Y_{n})-g(Y_{n})\big)^2
          \le \eps^2 r_n v_n, \quad \forall\, 1\le k\le N_\eps,
    $$
    then (D3) and (D4) can be replaced with the simpler
    condition
    $$ \int_0^\delta \sqrt{\log N_\eps}\, d\eps<\infty $$
    for some $\delta>0$ (cf.\ Theorem 2.11.9 of van der Vaart and
    Wellner, 1996).
    \item If $F(Y_{n})$ satisfies the Lindeberg condition (D2'),
    then (C2) and (D2) are satisfied.
    In particular, this holds if $nv_n\to\infty$ and
    \begin{equation} \label{ljapcond}
      E^* F(Y_{n})^{2+\delta} = O(r_nv_n) \quad \text{for some }
      \delta>0.
    \end{equation}
\item
    Thus, if (B1), (B2), (C3), (D1), (D3) and (D4) hold with a bounded envelope function
    $F$, then the empirical processes $Z_n$ converge to a centered
    Gaussian process with covariance function $c$.
  \end{enumerate}
  \hfill $\Box$
\end{remark}

\subsection{Asymptotic Equicontinuity}
  Like tightness, the asymptotic equicontinuity of $Z_n$ w.r.t.\
  $\rho$, i.e.
  $$ \forall\, \eps,\eta>0\; \exists \, \delta>0:  \quad
  \limsup_{n\to\infty} P^*\Big\{ \sup_{f,g\in\FF, \rho(f,g)<\delta}
  | Z_n(f)-Z_n(g)| >\eps\Big\}<\eta,
  $$
  is necessary and sufficient for the convergence of $Z_n$, provided
  all fidis of $Z_n$ converge.

To prove tightness we need a technical measurability condition,
Condition (D5) below, and, crucially, suitable bounds (D6) or (D6')
on the rate of increase of covering numbers. The condition (D5) in
particular is satisfied if the processes $(f(Y_{n}))_{f\in\FF}$ are
separable. The condition (D6) is stated in terms of a ``random
entropy'', while (D6'), which implies (D6), is phrased in terms of
uniform entropy. To state the assumptions, we need the following
definitions:

 For a given semi-metric $d$ on $\FF$, the (random) {\em covering number}
   $N(\eps,\FF,d)$ is the minimum number of balls with radius $\eps$
   w.r.t.\ $d$ needed to cover $\FF$. The condition (D6) bounds the rate of increase of
   $N(\eps,\FF,d_n)$ as $\eps$ tends to 0 for the random semi-metric
   $$ d_n(f,g) := \Big( \frac 1{nv_n} \sum_{j=1}^{m_n} \big(
   f(Y_{n,j}^*)-g(Y_{n,j}^*)\big)^2\Big)^{1/2},
   $$
   that is the $L_2$-semi-metric w.r.t.\ to empirical measure
   $(nv_n)^{-1}\sum_{j=1}^{m_n} \eps_{Y_{n,j}^*}$,
   where $Y_{n,j}^*$, $1\le j\le m_n$, are i.i.d.\ copies of $Y_{n,1}$.
    In    (D6') we instead use the supremum of all covering numbers $N(\eps,\FF,d_Q)$
    where
     $d_Q(f,g):= \big(\int (f-g)^2\, dQ\big)^{1/2}$ and $Q$
   ranges over the set of discrete probability measures
   $\mathcal{Q}$. With this notation, the conditions are as follows.
\\[1ex]
   {\bf (D5)}\;\;\;  \parbox[t]{11.4cm}{For all $\delta>0, n\in\N, (e_i)_{1\le i\le
   \floor{m_n/2}} \in \{-1,0,1\}^{\floor{m_n/2}}$ and $k\in\{1,2\}$ the
   map $\sup_{f,g\in\FF, \rho(f,g)<\delta}$
   $\sum_{j=1}^{\floor{m_n/2}} e_j\big(
   f(Y_{n,j}^*)-g(Y_{n,j}^*)\big)^k$
      is measurable.}\\[1ex]
\\[1ex]
   {\bf (D6)} \begin{equation*}
    \lim_{\delta\downarrow 0} \limsup_{n\to\infty}P^*\Big\{ \int_0^\delta
    \sqrt{ \log N(\eps,\FF,d_n)}\, d\eps > \tau \Big\} = 0, \quad
    \forall \tau>0.
  \end{equation*}
\\[1ex]
   {\bf (D6')}\;\;\;  \parbox[t]{11.2cm}{The envelope function $F$ is
   measurable with $E( F(Y_{n})^2)=O(r_nv_n)$ and  \begin{equation*}
    \int_0^1 \sup_{Q\in\mathcal{Q}} \sqrt{ \log
    N(\eps{\textstyle(\int F^2dQ)^{1/2}},\FF,d_Q)}\, d\eps < \infty.
  \end{equation*}}

\begin{theorem}  \label{condrem2}
  Suppose the basic assumptions (B1) and (B2) hold and that (D1), (D2'), (D3) and (D5)
  are satisfied. Then if also (D6) or (more restrictively, (D6')) holds, it follows that
    $Z_n$ is asymptotically equicontinuous. Further,  if in addition
    the finite-dimensional distributions converge (which in
    particular holds if (C1) and (C3) also are satisfied), then $Z_n$ converges to a
  Gaussian process with covariance function $c$ .
  \end{theorem}

 \begin{remark}\rm  \label{rem:vc} In view of (D6') one can apply the powerful Vapnik-\v{C}ervo\-nen\-kis
  theory to verify asymptotic equicontinuity. In particular,
  (D6') is satisfied if $\FF$ is a co-called
  VC-class or, more generally, a VC-hull class. We refer to Section
  2.6 of van der Vaart and Wellner (1996) for an outline of the most
  important uniform bounds on covering numbers $N(\eps(\int F^2dQ)^{1/2},\FF,d_Q)$.
     \hfill $\Box$
\end{remark}
\hspace*{4cm}

\section{Generalized tail array sums}
\label{sect:tailarrays}

Generalizing the tail empirical process $e_n(x)$ (for some fixed
$x\ge 0$), Rootz\'{e}n et al.\ (1990) considered so-called tail array
sums
\begin{equation} \label{tailarray}
 \sum_{i=1}^n   \phi(X_{n,i})
\end{equation}
for functions $\phi:\R\to\R$ satisfying $\phi(0)=0$ and $X_{n,i}$
defined by \eqref{standexceed}; see also Leadbetter and Rootz\'{e}n
(1993), Leadbetter (1995) and Rootz\'{e}n et al.\ (1998).

Like the tail empirical process, these tail array sums do not allow
inference about the extremal dependence structure, as the summands
$\phi(X_{n,i})$ depend on just one observation. However, if
$X_{n,i}$ denotes the vector of $d$ consecutive standardized
excesses, i.e.
\begin{equation} \label{standexcessvect}
  X_{n,i} := \bigg(\Big(
  \frac{X_i-u_n}{a_n}\Big)_+,\Big(\frac{X_{i+1}-u_n}{a_n}\Big)_+,\ldots,
  \Big(\frac{X_{i+d-1}-u_n}{a_n}\Big)_+\bigg),
\end{equation}
 then the statistic \eqref{tailarray} with
$\phi:(E,\mathcal{B}(E))\to(\R,\B)$ (and $E=\R^d$) contains
information on the extremal dependence structure.

Therefore, in the general setting of a row-wise stationary
triangular array $(X_{n,i})_{n\in\N,1\le i\le n}$ used in Section 2,
the {\em generalized (standardized) tail array sum} ({\em tail array
sum} for short) given by a measurable function
$\phi:(E,\mathcal{B}(E))\to(\R,\B)$ with $\phi(0)=0$ is defined as
\begin{equation}  \label{gentailarray}
  \tilde Z_n(\phi) := \frac 1{\sqrt{nv_n}}  \sum_{i=1}^n
  \big(\phi(X_{n,i})-E\phi(X_{n,i})\big).
\end{equation}
The tail array sum \eqref{gentailarray} can be obtained as the
empirical process $Z_n$ evaluated at the cluster functional
$$ g_\phi: E_\cup\to\R, \quad x=(x_1,\ldots,x_k)\mapsto
\sum_{i=1}^k \phi(x_i)
$$
if $n$ is a multiple of $r_n$. In general,  $\tilde
Z_n(\phi)-Z_n(g_\phi)=(nv_n)^{-1/2} \sum_{i=r_n m_n+1}^n$
$\big(\phi(X_{n,i})-E\phi(X_{n,i})\big)$, which is asymptotically
negligible under weak conditions specified in Corollary
\ref{gentailcorol} below.

For the remainder of this section, we assume that a family $\Phi$ of
functions $\phi$ of the above type is given, and assume it is
totally bounded w.r.t.\ a semi-metric $\rho_\Phi$ and has a finite
envelope function $\phi_{\max}:= \sup_{\phi\in\Phi}|\phi|$.
\begin{example}  \label{tailedfexam}  \rm ({\em Multivariate tail
empirical processes})
  If $X_{n,i}$ is defined as in \eqref{standexcessvect} and
  $\Phi:=\{ \boldsymbol{1}_{(x,\infty)} \mid x\in[0,\infty)^d\}$, then
  $\big(Z_n(g_\phi)\big)_{\phi\in\Phi}$ is the (reparametrized)
  multivariate tail empirical process. In particular, if $d=1$, then
  $\big(Z_n(g_\phi)\big)_{\phi\in\Phi}$  is a reparametrization
  of the tail empirical process  $e_n$ discussed in the introduction.

  For simplicity, we will assume that the $X_i$ are uniformly
  distributed; the general case can be easily obtained by a marginal
  quantile transformation (cf.\ Rootz\'{e}n (2009) for details). Then
  one chooses  $a_n=1-u_n=v_n$ for a sequence of thresholds $u_n$ tending to 1, so that the conditional
  distribution of the standardized excesses $X_{n,i}=(X_{i}-u_n)/a_n$ given that
  they are strictly positive is also uniform. Thus it suffices to
  consider $\Phi:=\{ \boldsymbol{1}_{(x,1]} \mid x\in[0,1]^d\}$ with envelope
  function $\phi_{\max}=\boldsymbol{1}_{(0,1]^d}$ and metric $\rho_\Phi( \boldsymbol{1}_{(x,1]},
  \boldsymbol{1}_{(y,1]}):= \max_{1\le l\le d} |x_l-y_l|$, $x,y\in[0,1]^d$.
  \hfill $\Box$
\end{example}

\begin{example}  \label{upcrossingexam}  \rm ({\em Upcrossings})
  If one is interested in upcrossings of a univariate time series over intervals
  $[x,y]$, then one may define $X_{n,i}$ as in Example
  \ref{tailedfexam} with $d=2$ and consider $\Phi:=\{
  \boldsymbol{1}_{[0,x)\times(y,1]} \mid x,y\in[0,1], x\le y\}$ with envelope
  function $\Ind{(x,y)\in [0,1]^2 \mid x < y}$.
  \hfill $\Box$
\end{example}

\begin{example}  \label{totalclaimsexam}  \rm ({\em Compound insurance
claim}) If $X_i$ denotes the $i$th claim of an insurance portfolio
with deductible $u_n+a_n t$ and $X_{n,i}$ as in \eqref{standexceed},
then $\phi_t:\R\to[0,\infty)$ given by
$\phi_t(x)=(x-t)\boldsymbol{1}_{(t,\infty)}(x)$ is  the standardized
total claimed amount. Thus the empirical process $\big(
 Z_n(g_{\phi_t})\big)_{t\ge 0}$ corresponding to
 $\Phi:=\{(x-t)\boldsymbol{1}_{(t,\infty)}(x) \mid t \geq 0\}$ describes the influence of the
 deductible on the random amount the insurance has to pay.
 \hfill $\Box$
\end{example}

\begin{example} \rm ({\em Bootstrapping the Hill estimator)}\label{ex:hill}
  A stationary  time
  series $(X_i)_{i\in\N}$ has extreme value index $\gamma>0$ if its marginal survival function $\bar F$ is regularly
  varying with index $-1/\gamma$, i.e.\ if  $\lim_{t\to\infty}\bar
  F(tx)/\bar F(t)=x^{-1/\gamma}$. Let $X_{n,i} := X_i/u_n 1_{\{X_i>u_n\}}$, $\phi_1(x)=\log (x)
  1_{\{x>1\}}$ and $\phi_2(x)=1_{\{x>1\}}$ so that $E \phi_2(X_{n,1})=v_n$ and $\gamma_n = E
  \phi_1(X_{n,1})/ E \phi_2(X_{n,1}) = E
  \phi_1(X_{n,1})/v_n=E(\log(X_1/u_n)|X_1>u_n)\to\gamma$ (cf.\ de
  Haan and Ferreira (2006), Theorem 1.2.1 and Remark 1.2.3). Then the Hill estimator $\hat{\gamma}_n$ of $\gamma$
  may be written as
\begin{equation} \label{eq:hill}
  \hat\gamma_n := \frac{\sum_{i=1}^n \log (X_i/u_n)
  1_{\{X_i>u_n\}}}{\sum_{i=1}^n   1_{\{X_i>u_n\}}}
  =\frac{\gamma_n + \tilde Z_n(\phi_1)/\sqrt{nv_n}}{
  1 + \tilde Z_n(\phi_2)/\sqrt{nv_n}}.
  \end{equation}
  Write  $g_k:=g_{\phi_k}$, $k\in\{1,2\}$, and suppose we  draw independent blocks $Y_i^{(n)}$
   from the empirical distribution of $Y_{n,i}$,
  $1\le i\le m_n$.
Then a bootstrap version of the Hill
  estimator is obtained as
  $$ \hat\gamma_n^* := \frac{\sum_{i=1}^{m_n} g_1(Y_i^{(n)})}{\sum_{i=1}^{m_n}
  g_2(Y_i^{(n)})}
  $$
\hfill $\Box$
\end{example}

\begin{example}  \label{densityexam}  \rm ({\em Kernel density
estimators})
  In this simple example we de\-mon\-strate that applications of the theory presented in
  Section 2 are not restricted to extreme value theory.
Further examples may be obtained from the literature on ``local
empirical processes''. For the analysis of such processes for
i.i.d.\ data we refer to Einmahl (1997), Gin\'e et al.\ (2003)  and
Gin\'{e} and Mason (2008) and to the lists of references in these
papers.

  Suppose that $(X_i)_{i\in\N}$ is a univariate
  stationary time series whose marginal df $H$ has a Lebesgue
  density $h$. Kernel estimators of the type
  $$ \hat h_n(x_0) := \frac 1{nb_n} \sum_{i=1}^n K\Big(
  \frac{X_i-x_0}{b_n} \Big)
  $$
  are probably the most widely used nonparametric estimators for
  $h(x_0)$ ($x_0\in\R$). Here $K$ denotes a suitable kernel, e.g.\
  a probability density with support $[-1,1]$, and $(b_n)_{n\in\N}$
  is a sequence of bandwidths tending to 0. Let
  $$ X_{n,i} := \Big( 2+\frac{X_i-x_0}{b_n}\Big)
  \boldsymbol{1}_{[x_0-b_n,x_0+b_n]}(X_i), \quad 1\le i\le n,
  $$
  where the constant 2 has been inserted to ensure $X_{n,i}>0$ for
  $X_i \in [x_0 -b_n, x_0+b_n]$. Let $\hat H_n$ be the corresponding empirical df. Then integration
  by parts yields
  \begin{eqnarray*}
    \hat h_n(x_0) & = & \frac 1{b_n} \int K(y-2)\, \hat H_n(dy) \\
      & = & \frac 1{b_n} \int \big(1 -\hat H_n(y+2)\big) \,K(dy)\\
      & = & \frac 1{nb_n} \int  \sum_{i=1}^n
      \boldsymbol{1}_{(y+2,\infty)}(X_{n,i})\, K(dy),
  \end{eqnarray*}
  provided that $K$ has bounded variation. Hence, for $\bar
  Z_n(y)=\tilde Z_n(\boldsymbol{1}_{(y+2,\infty)})$, $y\in[-1,1]$, and
  $n=r_nm_n$, we have that
  $$ \int \bar Z_n(y)\, K(dy) = \sqrt{\frac{n}{v_n}} b_n \big(
  \hat h_n(x_0)-E \hat h_n(x_0)\big),
  $$
  where $\sqrt{n/v_n} b_n\sim \sqrt{n/(2h(x_0) b_n)} b_n =
  \sqrt{nb_n/(2h(x_0))}$ as $n\to\infty$, if $h$ is continuous and positive at
  $x_0$. Thus one obtains the asymptotic normality of $\hat
  h_n(x_0)$ from the convergence of $\bar{Z}_n$ (or $\tilde Z_n$) towards a Gaussian
  process. Indeed, this way it is not difficult to derive normal
  approximations for $\hat h_n$ uniformly over families of kernels
  with compact support.
  \hfill $\Box$
\end{example}

To obtain conditions for weak convergence of tail array sums, we
first focus on families $\Phi$ such that the envelope function
$\phi_{\max}$ is bounded, which is true in the Examples
\ref{tailedfexam}, \ref{upcrossingexam} and \ref{densityexam}, but
not in Example \ref{totalclaimsexam} (unless the support of
$X_{n,i}$ is uniformly bounded). We  let $\FF:= \{g_\phi \mid
  \phi\in\Phi\}$ be equipped with the semi-metric $\rho(g_\phi,
  g_\psi)=\rho_\Phi(\phi,\psi)$.

\begin{corollary}    \label{gentailcorol}
  Suppose that $\phi_{\max}=\sup_{\phi\in\Phi}|\phi|$ is bounded and
  measurable, that $\Phi$ is totally bounded w.r.t.\ $\rho_\Phi$, that (B1) and (B2) hold,
   and that $r_n = o(\sqrt{nv_n})$. Further assume that
  \begin{equation}
        E\Big( \sum_{i=1}^{r_n} \Ind{X_{n,i}\ne 0}\Big)^2  = O(r_n
        v_n).
     \label{clustersizebound}
  \end{equation}
Then the conditions (C1),  (D1), and (D2') hold, and thus
    also (C2) and (D2) are satisfied.
 Moreover,
  \begin{equation}   \label{gentailarrayapprox}
     \sup_{\phi\in\Phi} \big| \tilde Z_n(\phi)-Z_n(g_\phi)\big| \to
     0 \quad \text{in outer probability.}
  \end{equation}
  If, in addition, (C3) and holds and one of the following two sets
  of  conditions
  \begin{itemize}
   \item[(i)]
  (D4) with a partition of $\FF$ independent of
    $n$,   or
    \item[(ii)]
    (D3), (D5),  and (D6)
  \end{itemize}
   are satisfied, then $\big(\tilde Z_n(\phi)\big)_{\phi\in\Phi}$, and the empirical
  processes $\big(Z_n(g_\phi)\big)_{\phi\in\Phi}$ of cluster
  functionals,   converge weakly to a Gaussian process with covariance
  function $c$.
\end{corollary}

\begin{remark}\rm  \label{covrem}
(i) It is possible to replace (C3) in the corollary by more basic
assumptions. Specifically, assume that the cluster lengths $L(Y_n)$
satisfy
\begin{equation} \label{eq:ltightness}
\lim_{k \to \infty} \limsup_{n \to \infty} \frac{1}{r_nv_n}
P\{L(Y_n)> k\}  = 0,
\end{equation}
 that there exist functions
$d_j:\Phi^2\to\R$ such that, for $k \in\N$ and $\phi,\psi\in\Phi$,
\begin{equation} \label{eq:phicovconv0}
\frac{1}{v_n}E\big(\phi(X_{n,1})\psi(X_{n,k})\big) \to d_{k-1}(\phi,
\psi),
  \; \text{as } n \to \infty,
\end{equation}
and that
\begin{equation} \label{eq:twoplusdeltabound}
E\Big( \sum_{i=1}^{r_n} \Ind{X_{n,i}\ne 0}\Big)^{2+\delta}  = O(r_n v_n),
\end{equation}
for some $\delta > 0$.  Then (C3'), and hence, by Corollary
\ref{covcondcorol}, also (C3) hold with
\begin{equation} \label{eq:dsum}
c(g_\phi, g_\psi) = d_0(\phi, \psi) + \sum_{i=1}^\infty \big( d_i(\phi, \psi) + d_i(\psi, \phi) \big).
\end{equation}
Moreover, if the following condition is met
\smallskip

\parbox[t]{1cm}{$\mathbf{(\widetilde{B3})}$} \parbox[t]{11.5cm}{For all $n\in\N$ and all $1\le i\le r_n$ there exists $s_n(i)\ge P(X_{n,i+1}\ne 0\mid X_{n,1}\ne 0)$ such that $s_\infty(i):=\lim_{n\to\infty} s_n(i)$ exists and $\lim_{n\to\infty} \sum_{i=1}^{r_n} s_n(i) =\sum_{i=1}^\infty s_\infty(i)<\infty$}
\smallskip

then \eqref{eq:ltightness} holds, and if, in addition, (B1) and  \eqref{eq:phicovconv0} are satisfied, then (C3) follows.
The proof is given in Section \ref{proofs}.

\noindent (ii) Suppose that the following simpler version of (C3'')
is satisfied, viz.\ that there exists a sequence $(W_i)_{i\in\N}$ of
$E$-valued random variables such that, for all $k\in\N$,
$P^{(X_{n,1},X_{n,k})\mid X_{n,1}\ne 0} \to P^{(W_1,W_k)}$ weakly,
with $P\{W_k\in D_\phi\setminus\{0\}\}=0$ for all $\phi\in\Phi$,
$k\in\N$, where $D_\phi$ is the discontinuity set of $\phi$. Then,
in view of Lemma \ref{lemma:clusterconv}, Remark \ref{condC3ppgen}
and  the boundedness of $\phi$ and $\psi$,
\begin{eqnarray*}
  \lefteqn{\frac{1}{v_n}E\phi(X_{n,1})\psi(X_{n,k}) =
E\big(\phi(X_{n,1})\psi(X_{n,k}) | \; X_{n,1} \neq 0 \big)}\\
  & \to & E\phi(W_1)\psi(W_k) =: d_{k-1}(\phi, \psi)
\end{eqnarray*}
so that equation \eqref{eq:phicovconv0} holds. \hfill $\Box$
\end{remark}

\begin{example} \label{tailedfexam2} \rm ({\em Multivariate tail empirical processes}, ctd)
In this example we give a set of conditions for the convergence of
the multivariate tail  empirical process from Example
\ref{tailedfexam} for uniformly  distributed rv's $X_i$. We then
discuss how  the condition (C3) on convergence of covariances may be
checked in the present situation. Finally we show that the central
condition \eqref{momcond} may be weakened in the univariate case, to
Condition \eqref{momcondd1}. This improves earlier results in the
literature.

Thus,  we first show that if $r_n = o(\sqrt{nv_n})$, (B1), (B2) and
(C3) are satisfied, and there exist a constant
$K$ and a $\delta>0$ such that for all sufficiently large $n$
  \begin{eqnarray}
    \label{momcond}\lefteqn{E\Big(\sum_{i=1}^{r_n} \boldsymbol{1}_{(x,y]}
    \Big(\frac{X_i-u_n}{a_n}\Big)\Big)^2 \le K|\log(y-x)|^{-(1+\delta)} r_nv_n,}\\
     &  &
     \hspace*{4cm} \forall \, 0\le x<y\le 1,\;
    y-x\le 1/2,\nonumber
  \end{eqnarray}
then the multivariate tail empirical process
$$  \bigg( \frac 1{\sqrt{nv_n}} \sum_{i=1}^n \Big(
  \boldsymbol{1}_{(x,1]}(X_{n,i})-P(X_{n,i}\in(x,1])\Big)\bigg)_{x\in[0,1]^d}
  $$ converges weakly to a Gaussian process with covariance function
  $c$.

Clearly \eqref{momcond} implies \eqref{clustersizebound}. By
Corollary \ref{gentailcorol}, it is hence enough to show that
Condition (i) of the corollary is satisfied. Now, to each
  $\eps>0$, let $\eta=\eta_\eps:=
  \exp\big(-(K^{-1}d^{-3}\eps^2)^{-1/(1+\delta)}\big)$ and define
  sets
  \begin{eqnarray*}
   \Phi^\eps_{(i_1,\ldots,i_d)} & := & \Big\{\boldsymbol{1}_{\times_{l=1}^d
  (x_l,1]} \mid (i_l-1)\eta\le x_l\le \min(i_l\eta,1)\forall\, 1\le l\le
  d\Big\},\\
  & & \hspace*{3cm}  \quad i_1,\ldots,i_d\in\{1,\ldots,\ceil{1/\eta}\},
  \end{eqnarray*}
  such that $\bigcup_{i_1,\ldots,i_d\in\{1,\ldots,\ceil{1/\eta}\}}
  \Phi^\eps_{(i_1,\ldots,i_d)}=\Phi$. Since, by (B1) and
  \eqref{momcond},
  \begin{eqnarray*}
  \lefteqn{E \sup_{\phi,\psi\in \Phi^\eps_{(i_1,\ldots,i_d)}}
  |g_\phi(Y_{n})- g_\psi(Y_{n})|^2}\\
   & = &
  E \Big( \sum_{i=1}^{r_n} \boldsymbol{1}_{\times_{l=1}^d
  ((i_l-1)\eta,1]\setminus \times_{l=1}^d
  (i_l\eta,1]} (X_{n,i})\Big)^2\\
  & \le & E \Big( \sum_{i=1}^{r_n} \sum_{l=1}^d \boldsymbol{1}_{((i_l-1)\eta,i_l\eta]}
  \Big(\frac{X_{i+l-1}-u_n}{a_n}\Big)\Big)^2\\
  & \le & d^2 E \max_{1\le l\le d}\Big(
    \sum_{i=1}^{r_n} \boldsymbol{1}_{((i_l-1)\eta,i_l\eta]}\Big(\frac{X_{i+l-1}-u_n}{a_n}\Big)\Big)^2\\
  & \le & d^3 K |\log \eta|^{-(1+\delta)} r_n v_n \\
  & = & \eps^2 r_n v_n,
  \end{eqnarray*}
  it follows that
  $$ \log N_{[\cdot]}(\eps,\FF,L_2^n) \le
  \log\big(\ceil{1/\eta}^d\big) = O(\eps^{-2/(1+\delta)})
  $$
  as $\eps\downarrow 0$. Hence the condition (D4) on entropy with
  bracketing holds with a partition independent of $n$, as required to prove the claim.

 The convergence  (C3) of covariance functions which was used above may sometimes be replaced
by simpler conditions. Specifically, Remark \ref{covrem}  gives
sufficient conditions for (C3) to hold, for general $d\in\N$. Assume
e.g.\ that
  all bivariate distributions $(X_1, X_m)$ belong to the domain of
  attraction of some bivariate extreme value distribution.
  Then, since the limiting random variables $W_i$ are continuous on $(0,\infty)$,
  the assumptions of Remark \ref{covrem} (ii) are satisfied, and hence \eqref{eq:phicovconv0} holds
  (cf.\ Segers, 2003, Theorem 2).   Further, Condition \eqref{eq:twoplusdeltabound} holds
  if and only if for some $\delta>0$
  \begin{equation} \label{ljapcond3}
     E\Big( \sum_{i=1}^{r_n} \boldsymbol{1}_{(u_n,1]} (X_i)\Big)^{2+\delta}
  = O(r_nv_n).
  \end{equation}

 For the case $d=1$, the condition \eqref{momcond} can be
  weakened, to the requirement that
  \begin{equation}   \label{momcondd1}
    E\Big(\sum_{i=1}^{r_n} \boldsymbol{1}_{(x,y]}
    \Big(\frac{X_i-u_n}{a_n}\Big)\Big)^2 \le
    h(y-x) r_nv_n \quad \forall \, 0\le x<y\le 1,
  \end{equation}
  for some function $h:(0,\infty)\to(0,\infty)$ satisfying
  $\lim_{t\downarrow 0} h(t)=0$. To see this, note that the
  functions $\phi_x=\boldsymbol{1}_{(x,1]}$, $x\in[0,1]$, are linearly ordered,
  and hence so are the corresponding cluster functionals $g_{\phi_x}$,
  $x\in[0,1]$. Hence $\FF=\{g_{\phi_x} \mid  x\in[0,1]\}$ is
   a VC class of functions (van der Vaart and
  Wellner (1996), Section 2.6). Thus,
  according to Remark \ref{rem:vc}, (D6') (and hence also (D6)) is satisfied. The measurability
  condition (D5) holds, since all processes
  occurring in this setting are separable. Moreover, (D3)
  is satisfied for the metric
  $\rho(g_{\phi_x},g_{\phi_y}):=|y-x|$:
  \begin{eqnarray*}
    \lefteqn{\limsup_{n\to\infty} \frac 1{r_nv_n} \sup_{x,y\in[0,1],
    |y-x|<
    \delta} E\big(g_{\phi_x}(Y_{n})-g_{\phi_y}(Y_{n})\big)^2}\\
    & = & \limsup_{n\to\infty} \frac 1{r_nv_n} \sup_{x,y\in[0,1],
    |y-x|<
    \delta} E\Big( \sum_{i=1}^{r_n} \boldsymbol{1}_{(x,y]} (X_{n,i})\Big)^2\\
    & \le & \sup_{0<t\le\delta} h(t)\\
    & \to & 0
  \end{eqnarray*}
  as $\delta\downarrow 0$ by \eqref{momcondd1}, so that version (ii)
of Corollary \ref{gentailcorol} applies. This proves the claim that
\eqref{momcond} may be weakened to \eqref{momcondd1} in the
univariate case.

If we could assume that $\{X_i; 1\leq i \leq n \}$ could be split up
in into consecutive {\em independent} blocks of length $r_n$ then
\eqref{momcondd1} would  be seen to be the same as to assume that
$E\big(Z_n(g_{\phi_y}) - Z_n(g_{\phi_x})\big)^2 \leq h(|y-x|)$, for
some $h$ with properties as above. This is the same as to assume
that $Z_n$ is uniformly mean square continuous. However, in the
proofs in Section \ref{proofs} we use mixing to translate to cases
where this independence assumption in fact can be made, and
accordingly \eqref{momcondd1} seems quite minimal. In fact, in view
of the counterexamples in Hahn (1977) it may even be surprising that
this condition is sufficient.

  Rootz\'{e}n (1995, 2009) proved convergence of the univariate tail empirical
  process $e_n$  using a more restrictive version of \eqref{momcond} and the stronger
  condition that $r_n=o((nv_n)^{1/2-\eps})$ for some $\eps>0$. In Drees (2000) Rootz\'en's conditions
  were  slightly weakened to the requirements
  that $r_n=o((nv_n)^{1/2}\log^{-2}(nv_n))$ and that
  \begin{equation}   \label{rootzencond}
    E\Big(\sum_{i=1}^{r_n} \boldsymbol{1}_{(x,y]}
    \Big(\frac{X_i-u_n}{a_n}\Big)\Big)^2 \le
    K(y-x) r_nv_n \quad \forall \, 0\le x<y\le 1,
  \end{equation}
  instead of \eqref{momcond}. Condition
  \eqref{rootzencond} is much more restrictive than \eqref{momcond}
  for small $y-x$. In many specific time  series models,
  it was condition  \eqref{rootzencond} (for small $y-x$) that turned out to be most difficult
  to verify; see e.g.\ the discussion of the solutions of a
  stochastic recurrence equation in Drees (2000), Section 4.
  Therefore, it might be useful that the bound in \eqref{momcond}
  converges to 0 much more slowly as $y-x$ tends to 0.
  \hfill $\Box$
\end{example}

It is possible to deal with  Examples \ref{upcrossingexam} and
\ref{densityexam} in a similar fashion.

As already mentioned, Example \ref{totalclaimsexam} does not fit
into the framework of Corollary \ref{gentailcorol} if the underlying
df belongs to the domain of attraction of an extreme value
distribution with non-negative extreme value index, because then the
support is not bounded. In that case, Condition
\eqref{clustersizebound} must be strengthened.
\begin{corollary} \label{gentailcorol2}
  In the setting of Corollary \ref{gentailcorol} the assertions
  remain true if $\phi_{\max}$ is measurable but not necessarily
  bounded, provided \eqref{clustersizebound} is replaced with
  \begin{equation}  \label{clustersumbound}
    E\Big( \sum_{i=1}^{r_n} \phi_{\max}(X_{n,i})\Big)^{2+\delta} =
    O(r_nv_n) \quad \mbox{for some } \delta>0.
  \end{equation}
\end{corollary}

\begin{example} \rm ({\em Compound insurance
claim}, ctd) In the setting of Example \ref{totalclaimsexam},
  uniform convergence of the empirical process of cluster
  functionals can be expected only if the deductible $t$ is
  restricted to some bounded set. Therefore, we consider the set $\Phi_T:=\{
  \phi_t \mid t\in[0,T]\}$ for an arbitrary  $T\in(0,\infty)$ This
  set  is totally bounded w.r.t.\ the metric
  $d_\Phi(\phi_s,\phi_t):=|s-t|$. The envelope function is
  $\phi_{\max}(x)=\phi_0(x)=x_+$.

Suppose the conditions (B1), (B2), (C3), \eqref{clustersizebound},
and
\begin{equation} \label{momcond3}
 E\Big(\sum_{i=1}^{r_n} X_{n,i}\Big)^{2+\delta}
=O(r_nv_n)
\end{equation}
  for some $\delta>0,$
  are satisfied. Then the empirical process
  $(Z_n(g_{\phi_t}))_{0\le t\le T}$ converges weakly to a Gaussian process.

  To see this, first observe that the functions $\phi_t$ are monotonically decreasing
  in $t$. Hence $\Phi_T$ is a VC class of functions, so that
  (D6) holds (see Remark \ref{condrem2}).
  Since all sample paths are continuous, the measurability condition
  (D5)  trivially holds.

  To prove (D3) check that
  \begin{eqnarray*}
     \lefteqn{\sup_{0\le s\le t\le T, |t-s|<\delta}  \frac
    1{r_nv_n} E\Big(\sum_{i=1}^{r_n}\big( (X_{n,i}-s)_+ -
    (X_{n,i}-t)_+\big)\Big)^2} \\
    & \le & \sup_{0\le s\le t\le T, |t-s|<\delta} \frac    1{r_nv_n} E\Big(\sum_{i=1}^{r_n}
    (t-s)\boldsymbol{1}_{(s,\infty)}(X_{n,i})\Big)^2\\
    & \le & \delta^2 \frac
    1{r_nv_n} E\Big( \sum_{i=1}^{r_n}
    \boldsymbol{1}_{(0,\infty)}(X_{n,i})\Big)^2.
  \end{eqnarray*}
  By \eqref{clustersizebound}, the $\limsup$ of the right-hand side (as $n$
  tends to $\infty$) is bounded by a multiple of $\delta^2$, which
  yields (D3). Further, \eqref{momcond3} is just a reformulation of
  \eqref{clustersumbound} to the present setting. Hence all the
  conditions of Corollary \ref{gentailcorol2} have been verified,
  and thus the result follows.

  By Corollary \ref{covcondcorol}, the condition (C3) in turn follows if, in
  addition, one assumes that all finite-dimensional
  marginal distributions of the time series $(X_i)_{i\in\N}$ belong to the domain of attraction of some
  extreme value distributions and that the normalizing constants $u_n$ and $a_n$ are
  chosen accordingly. Then  (C3.1'') holds  (cf.\ Segers, 2003, Theorem
  2), and
  (C3.2'') also follows, from \eqref{clustersumbound}  and Lemma \ref{simplerfidicond}
  (vi).  \hfill $\Box$
\end{example}

\begin{example} \rm ({\em Bootstrapping the Hill estimator}, ctd)
Continuing Example \ref{ex:hill} we now sketch proofs of asymptotic
normality of the Hill estimator and of consistency of the block
bootstrap.  Full process convergence may also be obtained and is
useful if e.g.\  $u_n$ is replaced by $k_n$-th largest order
statistic, for some suitable sequence $k_n$. We use asymptotic
normality to show consistency of the block bootstrap -- but the hope
is that the bootstrap has better small-sample properties than the
normal approximation with estimated variance.

For this we assume that (B1) and (B2) and, with the notation of Example \ref{ex:hill}, that for  $k,l \in  \{1,2\}$
  \begin{eqnarray}
     E\Big(\sum_{i=1}^{r_n}
     \phi_k(X_{n,i})\Big)^4 & = & O(r_nv_n)
     \label{eq:bootstrapcond1}\\
     \lim_{n\to\infty} \frac 1{r_n v_n} \sum_{i=1}^{r_n} \sum_{j=1}^{r_n}
     E\big(\phi_k(X_{n,i})\phi_l(X_{n,j})\big) & = & \sigma_{kl} \nonumber
  \end{eqnarray}

Then, in a similar way as in the proofs of
  Corollaries \ref{gentailcorol} and \ref{gentailcorol2}, it can be seen that
  $(\tilde Z_n(\phi_k))_{1\le k\le 2}$ converges to a centered
  normal distribution with covariance matrix $(\sigma_{kl})_{1\le
  k,l\le 2}$. It
  follows that
  \begin{equation} \label{eq:hillapprox}
   \hat\gamma_n = \gamma_n + (nv_n)^{-1/2}
  \big(\tilde Z_n(\phi_1)-\gamma \tilde Z_n(\phi_2)\big) +
  o_p\big((nv_n)^{-1/2}\big),
  \end{equation}
  and thus that
  \begin{equation} \label{Hillasymp}
   \sqrt{nv_n}(\hat\gamma_n-\gamma_n)\;\longrightarrow\;
  {\mathcal N}_{(0,\sigma_{11}+\gamma^2\sigma_{22}-2\gamma\sigma_{12})}, \; \; \; \mbox{in distribution.}
  \end{equation}

  Writing  $X^{(n)}:= (X_i)_{1\le i\le n}$ for the original data we next show that
\begin{equation} \label{eq:bootstrapapprox}
   \sup_{t\in\R}\Big| P\big(
  \sqrt{nv_n}(\hat\gamma_n^*-\hat\gamma_n)\le t\mid X^{(n)}\big)- P\big\{
  \sqrt{nv_n}(\hat\gamma_n-\gamma_n)\le t\big\}\Big|=o_P(1),
  \end{equation}
  i.e.\ consistency of the block bootstrap estimator.  With the notation from Example \ref{ex:hill},
  $$ \frac{E\big(g_1(Y_i^{(n)})|X^{(n)}\big)}{E\big(g_2(Y_i^{(n)})|X^{(n)}\big)}
  = \frac{m_n^{-1} \sum_{i=1}^{m_n} g_1(Y_{n,i})}{m_n^{-1} \sum_{i=1}^{m_n}
  g_2(Y_{n,i})} = \hat\gamma_n.
  $$
  From arguments as in the proof of Lemma \ref{iidcompprop} below (in particular
  \eqref{totvardist}), it follows that if  condition
  \eqref{eq:bootstrapcond1} holds then $Z_n(g_k g_l)=O_P(1)$. Hence, for $k,l\in\{1,2\}$,
  \begin{eqnarray*}
    \lefteqn{ \frac 1{r_nv_n} Cov\big( g_k(Y_1^{(n)})
    g_l(Y_1^{(n)})\mid X^{(n)} \big)} \\
    & = & \frac 1{r_nv_n} \Big( \frac 1{m_n} \sum_{i=1}^{m_n}  g_k(Y_{n,i})
    g_l(Y_{n,i}) - \frac 1{m_n} \sum_{i=1}^{m_n}  g_k(Y_{n,i})\cdot\frac 1{m_n} \sum_{i=1}^{m_n}
    g_l(Y_{n,i})\Big)\\
    & = & \frac 1{r_nv_n}Cov\big(g_k(Y_{n,1}),g_l(Y_{n,1})\big) - \frac 1{m_n} Z_n(g_k)Z_n(g_l)\\
    & & \hspace*{0.5cm} + \frac
    1{\sqrt{nv_n}} \big( Z_n(g_kg_l)-E(g_l(Y_{n,1}))Z_n(g_k)
    - E(g_k(Y_{n,1}))Z_n(g_l)\big)\\
    & \to & \sigma_{kl}
  \end{eqnarray*}
  in probability. Similarly as in \eqref{eq:hillapprox} we have  that
  \begin{eqnarray*} \label{bootstrapexpan}
  \hat\gamma_n^*
  & = & \hat\gamma_n +
  o_p\big((nv_n)^{-1}\big)\\
   & &  \hspace*{-0.4cm}+ (nv_n)^{-1}
  \sum_{i=1}^{m_n} \big(g_1(Y_i^{(n)})-\gamma g_2(Y_i^{(n)})-E(g_1(Y_i^{(n)})-\gamma g_2(Y_i^{(n)})|X^{(n)})\big).
  \end{eqnarray*}
  Moreover, one can conclude from \eqref{eq:bootstrapcond1} that
  \begin{eqnarray*}
   m_nE\Big( \Big(\frac{g_k(Y_1^{(n)})-
  E(g_k(Y_1^{(n)})|X^{(n)})}{\sqrt{nv_n}}\Big)^3\,\Big|\, X^{(n)}\Big) &
  = &
  O_P\big(m_n(nv_n)^{-3/2} r_nv_n\big)\\
  & = & O_P((nv_n)^{-1/2}),
  \end{eqnarray*}
  and thus  the Berry-Ess\'{e}en
  inequality yields
  \begin{eqnarray*}
    \lefteqn{\sup_{t\in\R} \Big| P\Big( (nv_n)^{-1/2} \sum_{i=1}^{m_n} \big(g_1(Y_i^{(n)})
  -\gamma g_2(Y_i^{(n)})-E(g_1(Y_i^{(n)})-\gamma
  g_2(Y_i^{(n)})|X^{(n)})\big)}\hspace*{2cm}\\
    & & \hspace*{0.8cm} \le t\mid X^{(n)}\Big)-
  \Phi\big((\sigma_{11}+\gamma^2\sigma_{22}-2\gamma\sigma_{12})^{-1/2}t\big)\Big|=o_P(1).
  \end{eqnarray*}
  In view of \eqref{Hillasymp}  this   proves \eqref{eq:bootstrapapprox}.
  \hfill $\Box$
\end{example}

\section{Indicator functionals}
\label{sect:indicators}

Another important class of cluster functionals are indicator
functions. Notice that by definition these indicator functions are
applied to whole clusters, while in the Examples \ref{tailedfexam},
\ref{upcrossingexam} and \ref{densityexam} above indicator functions of
single observations $X_{n,i}$ were summed up. For $C\subset E_\cup$
the indicator function $\boldsymbol{1}_C$ is a cluster functional if
and only if the set satisfies the following two conditions:
\begin{itemize}
  \item $x = (x_1, \ldots x_\ell) \in C \iff (0,x_1, \ldots x_\ell)\in C
  \iff (x_1, \ldots x_\ell,0)\in C$ for all $x\in   E_\cup$
  \item $0\not\in C$
\end{itemize}

In this section we study situations where the set of cluster
functionals is of the form\{$\FF=\{\boldsymbol{1}_C\mid C\in\CC\}$
for some family $\CC\subset 2^{E_\cup}$ of such sets.

\begin{example}  \rm \label{ex:survivalfunction} ({\em Joint survival
function of cluster values}) The conditional joint survival
  function of the first $k$ observations in a cluster core
  $Y_n^c$, given that the core has length greater than or equal
  to $k$, can be estimated by
  $$ \frac{\displaystyle \sum_{j=1}^{m_n} \boldsymbol{1}_{C_{t_1,\ldots,
  t_k}}(Y_{n,j})}{\displaystyle \sum_{j=1}^{m_n} \boldsymbol{1}_{C_{0,\ldots,
  0}}(Y_{n,j})}
  $$
  with
  $$ C_{t_1,\ldots, t_k} := \big\{ x\in E_\cup \mid \exists\, j:
  x_i=0\; \forall\, 1\le i\le j,\; x_{j+i}>t_i\; \forall\, 1\le i\le
  k\big\}.
  $$
  Obviously, a limit theorem for the empirical process
  $$ \tilde Z_n(t_1,\ldots, t_k) := Z_n(\boldsymbol{1}_{C_{t_1,\ldots,  t_k}}),
  \quad t_1,\ldots,t_k\in [0,1],
  $$
  is useful for the asymptotic analysis of the above estimator.
  \hfill $\Box$
\end{example}

\begin{example}  \rm \label{ex:orderstat} ({\em Order statistics of
cluster values}) Let
  $$ D_{t_1, \ldots t_k}:=  \bigcap_{j=1}^k E_{j,t_j}
  $$
  with
  $$ E_{j,t_j} := \big\{ (x_1,\ldots, x_m) \in E_\cup \mid m\in\N,
  \sum_{i=1}^m \boldsymbol{1}_{(t_j,1]}(x_i) \ge j\big\},
  $$
  i.e., $D_{t_1, \ldots t_k}$ contains all vectors of
  arbitrary length such that the $j$th largest value exceeds $t_j$
  for all $1\le j\le k$. Then the empirical process
  $\tilde{Z}_n(t_1, \ldots t_k) = Z_n(\boldsymbol{1}_{D_{t_1, \ldots t_k}})$
  describes the standardized joint empirical survival
  function of the $k$ largest order statistics of the cluster cores.
  \hfill $\Box$
\end{example}

Next we discuss the conditions imposed in Theorem \ref{condrem2} to
ensure convergence of the empirical processes considered in this
section.

 The conditions (D1) and (D2') are trivial, and condition
(C1) holds by Lemma \ref{simplerfidicond} (ii).

If $r_nv_n \to 0$ (which is a part of assumption (B1)), then (C3)
is equivalent to
\begin{equation}
   \frac 1{r_nv_n} P\{Y_{n,1}\in C\cap
  D\}
   \to  c(\boldsymbol{1}_C,\boldsymbol{1}_D),   \label{covcond2}
\end{equation}
since $\Cov(\boldsymbol{1}_C(Y_n),\boldsymbol{1}_D(Y_n)) =
P\{Y_{n}\in C\cap D\}-P\{Y_{n}\in
  C\}\cdot P\{Y_{n}\in D\}$ and since $P\{Y_{n}\in
  C\}\cdot P\{Y_{n}\in D\} = O((r_nv_n)^2)=o(r_nv_n)$.

Similarly, condition (D3) can be reformulated as
\begin{equation}   \label{ascontdiff}
  \lim_{\delta\downarrow 0} \limsup_{n\to\infty} \sup_{C,D\in\CC,
  \rho_\CC(C,D)<\delta} \frac 1{r_n v_n} P\{Y_{n}\in C\triangle
  D\}=0
\end{equation}
where $C\triangle  D=(C\setminus D) \cup (D\setminus C)$ denotes the
symmetric difference between $C$ and $D$ and $\rho_\CC$ is a
semi-metric on $\CC$ that induces a semi-metric $\rho$ on $\FF$ via
$\rho(\boldsymbol{1}_C,\boldsymbol{1}_D):=\rho_\CC(C,D)$.

If (C3'') holds, then
$$ \frac 1{r_n v_n} P\{Y_{n}\in C\triangle D\} \;
\longrightarrow\; P\{(W_i)_{i\ge 1} \in C\triangle D\} -
P\{(W_i)_{i\ge 2} \in C\triangle D\},
$$
where $(W_i)_{i\ge 1} \in C\triangle D$ is interpreted as
$(W_i)_{1\le i\le m} \in C\triangle D$ for some $m\ge m_W$, i.e.\
$W_i=0$ for all $i>m$. If the following continuity property holds
$$ \lim_{\delta\downarrow 0} \sup_{C,D\in\CC,\rho_\CC(C,D)<\delta} P\{(W_i)_{i\ge 1} \in C\triangle D\} -
P\{(W_i)_{i\ge 2} \in C\triangle D\} =0,
$$
then results by Fabian (1970) may help to conclude (D3). However, in
the  examples of this section we will verify (D3) in a more direct
way.

Finally, if $\CC$ is a VC-class, then condition (D6') is fulfilled
(cf.\ Remark \ref{rem:vc}).

The following result gives conditions for the convergence of the
empirical processes in  Examples \ref{ex:survivalfunction} and
\ref{ex:orderstat}. Here we assume  that the random variables
$X_{n,i}$ are $[0,1]$-valued so that is suffices to consider the
processes $\tilde Z_n$ with index set $[0,1]^k$. If the rv's
$X_{n,i}$ are standardized excesses defined in \eqref{standexceed}
(as we assume in the second part of the following corollary), then
this can be achieved by a simple quantile transformation (cf.\
Example \ref{tailedfexam}).

\begin{corollary} \label{cor:indicators}
 \begin{enumerate}
 \item
Let $ \tilde{Z}_n(t_1,\ldots, t_k) $ be as in Examples
\ref{ex:survivalfunction} or  \ref{ex:orderstat}, with $t_i \in [0,
1], i=1, \ldots k$, and suppose (B1), (B2), (B3), (C3.1''), and (D3)
hold with $\rho\big(\boldsymbol{1}_{C_{s_1,\ldots,
s_k}},\boldsymbol{1}_{C_{t_1,\ldots,  t_k}}\big):=\sum_{i=1}^k
|s_i-t_i|$ resp.\\ $\rho\big(\boldsymbol{1}_{D_{s_1,\ldots,
s_k}},\boldsymbol{1}_{D_{t_1,\ldots,  t_k}}\big):=\sum_{i=1}^k
|s_i-t_i|$. Then $\tilde{Z}_n$ converges to a continuous Gaussian
process. If $\tilde{Z}_n$ is as in Example
\ref{ex:survivalfunction}, then the covariance function of the
process is
 \begin{eqnarray} \label{ex:survivalcov}
    \lefteqn{\tilde c\big((s_1,\ldots,s_k),(t_1,\ldots, t_k)\big)=P\{(W_i)_{i\ge 1}
 \in C_{\max(s_1,t_1),\ldots, \max(s_k,t_k)}\}  }\hspace*{4cm} \nonumber \\
      &  & { } - P\{(W_i)_{i\ge 2}
 \in C_{\max(s_1,t_1),\ldots, \max(s_k,t_k)}\},
 \end{eqnarray}
 and if $\tilde{Z}_n$ is as in Example \ref{ex:orderstat}, then the covariance function of the process is
  \begin{eqnarray} \label{eq:ordercov}
    \lefteqn{\tilde c\big((s_1,\ldots,s_k),(t_1,\ldots, t_k)\big)= P\Big\{(W_i)_{i\ge 1}
 \in \bigcap_{j=1}^k E_{j,\max(s_j,t_j)}\Big\}}\hspace*{4cm} \nonumber \\
      &  &  - P\Big\{(W_i)_{i\ge 2}
 \in\bigcap_{j=1}^k E_{j,\max(s_j,t_j)}\Big\}.
\end{eqnarray}
 \item More specifically, assume that the rv's $X_{n,i}$ are standardized
 excesses of a uniformly distributed univariate stationary time
 series (as in Example \ref{tailedfexam}) and that all
 finite-dimensional marginal distributions belong to the domain of
 attraction of some extreme value distribution. Then the assertions
 of part (i) hold true if the conditions (B1), (B2) and (B3) are
 satisfied.
 \end{enumerate}
\end{corollary}

In Example \ref{ex:survivalfunction} we only considered the first
$k$ ``extremes'' in each cluster, where $k$ is a fixed number. Since
for most time series the cluster size is not bounded, the resulting
empirical process does not give a full picture of the stochastic
behavior of the clusters. To overcome this drawback, in the final
example we define and analyze an empirical process of cluster
functionals that takes {\em all} values of each cluster into
account. As the cluster length is random, this requires work with
a quite complex index set.

\begin{example}  \rm \label{ex:allvalues} ({\em Joint distribution of all cluster values})
Recalling the notation $L(x)$ for the length, say $j$, of the core
$x^c = (x_1^c, \ldots, x_j^c)$ of a vector $x$, we set
$$
C_{j, t_1, \ldots, t_j} := \{x \in E_\cup \mid L(x) = j, x_i^c \in
[0, t_i], \forall 1 \leq i \leq j \}.
$$
Then the empirical process $\tilde{Z}_n(j, t_1, \ldots, t_j) :=
Z_n(\boldsymbol{1}_{C_{j, t_1, \ldots, t_j}})$, $j\in\N$, $t_i\ge
0$, describes the joint distribution of  all the values in a
cluster.

Like in Corollary \ref{cor:indicators} (ii), for simplicity we focus
on the case that the clusters are based on standardized exceedances
$X_{n,i}$ of a uniformly distributed stationary time series
$(X_i)_{i\in\N}$, such that all finite-dimensional marginal
distributions belong to the domain of attraction of some extreme
value distribution. However, it is not difficult to generalize this
result to a slightly more general setting which is analog to the one
considered in Corollary \ref{cor:indicators} (i).

Suppose that  (B1), (B2), and (B3)  hold, and that
\begin{equation} \label{eq:lengthhalfbound}
E\big(L(Y_n)^{1 + \zeta} \mid Y_n \neq 0\big) = O_p(1), \;\;\;
\mbox{some} \;\; \zeta > 0.
\end{equation}
Then $\tilde{Z}_n$ converges weakly to a continuous Gaussian process with covariance function
\begin{eqnarray}
\lefteqn{c\big((j, s_1, \ldots, s_j),(k, t_1, \ldots,
t_k)\big)}\nonumber\\
\label{eq:allvaluescov} &= & \delta_{j,k}\Big(P\big\{L(W) = k, W_i \leq s_i \wedge t_i, \; \forall 1\leq i \leq k\big\} \\
 & &{ } - P\big\{L(W^{(2;\infty)})=k,\,
  \big((W^{(2;\infty)})^c\big)_i\le t_i,\, \forall 1\le i\le k \big\}
  \Big) \nonumber
\end{eqnarray}
where $\delta_{j,k}$ is one if $j=k$ and zero otherwise.

The proof of this uniform central limit theorem is given in Section
\ref{proofs}.

\hfill $\Box$
\end{example}

\section{Proofs} \label{proofs}

In this section we prove the results from Sections
\ref{sec:limittheorems}--\ref{sect:indicators}. We start with fidi
convergence, then consider asymptotic tightness and asymptotic
equicontinuity, and finally prove the corollaries from Sections
\ref{sect:tailarrays} and \ref{sect:indicators}.

The first step in the proof of fidi convergence is to use mixing to
bring the problem back to classical limit theory for iid variables.
Let $Y_{n,j}^*$ denote iid copies of the original blocks $Y_{n,j}$
(which are identically distributed, but are not assumed to be
independent -- and which in interesting cases typically are
dependent).

\begin{lemma}  \label{iidcompprop}
  Suppose (B1), (B2) and (C1) are satisfied. Then the fidis of $(Z_n(f))_{f\in\FF}$
  converge weakly if and only
  if the fidis of the sums of independent blocks
  $$ Z_n^*(f) := \frac 1{\sqrt{n v_n}} \sum_{j=1}^{m_n} \big( f(Y_{n,j}^*)-
E f(Y_{n,j}^*)\big), \quad f\in\FF,
$$
converge weakly. In this case the limit distributions are the same.
\end{lemma}

\begin{proof}
\rm\quad Let
\begin{eqnarray*}
      \Delta^*_{n,j}(f) & := &  f(Y_{n,j}^*) - f((Y_{n,j}^*)^{(r_n-l_n)}), \quad
   1\le j\le m_n,
\end{eqnarray*}
and let $ \Delta_{n,j}(f)$ be defined in the same way, but instead
based on the original (dependent) blocks, so that $
\Delta^*_{n,j}(f) \stackrel{d}{=} \Delta_{n,j}(f) \stackrel{d}{=}
\Delta_n(f)$ for each $j$,  with $\Delta_n(f)$ as in (C1). By
Theorem 1 in Petrov
 (1975), Section IX.1, applied to the iid random variables
 $X_{nk}:=(nv_n)^{-1/2}\Delta_{n,k}^*(f)$, condition (C1) implies that
 \begin{equation} \label{bbsbcond}
        \frac 1{\sqrt{n v_n}} \sum_{j=1}^{m_n} \big( \Delta_{n,j}^*(f)-
E \Delta_{n,j}^*(f)\big)=o_P(1), \quad \forall\, f\in\FF.
  \end{equation}

  We next prove  the analogous
  convergence for the dependent random variables, i.e. that
   \begin{equation} \label{bbsbdep}
        \frac 1{\sqrt{n v_n}} \sum_{j=1}^{m_n} \big( \Delta_{n,j}(f)-
E \Delta_{n,j}(f)\big)=o_P(1) \quad \forall\, f\in\FF.
  \end{equation}

 Using Theorem 1 in Petrov
 (1975), Section IX.1 again, it also follows from (C1) that the
  convergence analogous to \eqref{bbsbcond}  hold for the sums of the
  even numbered blocks
  \begin{eqnarray}
    \frac 1{\sqrt{n v_n}} \sum_{j=1}^{\floor{m_n/2}} \big( \Delta_{n,2j}^*(f)-
E \Delta_{n,2j}^*(f)\big)=o_P(1).  \label{bbsbeven}
\label{bbsbodd}
  \end{eqnarray}

  Since the even numbered blocks $Y_{n,j}$ are separated by $r_n$
  observations, a well-known inequality for the total variation
  distance (cf.\ Eberlein, 1984) between the joint distributions
  of dependent observations
  and independent copies  yields
  \begin{equation}  \label{totvardist}
    \big\| P^{(Y_{n,2j})_{1\le j\le
    \floor{m_n/2}}} - P^{(Y_{n,2j}^*)_{1\le j\le
    \floor{m_n/2}}} \big\|_{TV} \le \floor{m_n/2} \beta_{n,r_n}
    \to 0
  \end{equation}
  by (B2). Combining \eqref{bbsbeven} with
  \eqref{totvardist}, we arrive at
  $$  \frac 1{\sqrt{n v_n}} \sum_{j=1}^{\floor{m_n/2}} \big( \Delta_{n,2j}(f)-
E \Delta_{n,2j}(f)\big)=o_P(1).
  $$
  Together with the analogous convergence for the sum over the odd
  numbered blocks this proves \eqref{bbsbdep}.

  Thus the fidis of $Z_n$ converge if and only if the fidis of
  \begin{eqnarray*}
   \bar Z_n(f) & := & Z_n(f) - \frac 1{\sqrt{n v_n}} \sum_{j=1}^{m_n} \big( \Delta_{n,j}(f)-
E \Delta_{n,j}(f)\big)\\
    & = & \frac 1{\sqrt{n v_n}} \sum_{j=1}^{m_n} \big( f(
    Y_{n,j}^{(r_n-\ell_n)})-E f( Y_{n,j}^{(r_n-\ell_n)})\big), \quad f\in\FF,
  \end{eqnarray*}
  converge, and in this case the limiting distributions are the
  same. Similarly, by \eqref{bbsbcond} the corresponding assertion holds for the
  sums over the independent blocks, and then the lemma follows from the  inequality
  for the total variation distance, since it implies that
  $$ \big\| P^{(Y_{n,j}^{(r_n-\ell_n)})_{1\le j\le
    m_n}} - P^{((Y^*_{n,j})^{(r_n-\ell_n)})_{1\le j\le
    m_n}} \big\|_{TV} \le m_n \beta_{n,l_n}
    \to 0
  $$
 by (B2), since the shortened blocks $Y_{n,j}^{(r_n-\ell_n)}$ are separated
  by $l_n$ observations.
\end{proof}

\begin{proofof} {\sc Theorem \ref{fidisconvtheo}.}\rm\quad
  The assertion follows from Lemma \ref{iidcompprop} and
   and the multivariate central limit theorem
  for triangular arrays of row-wise independent random vectors applied to
  $(Z_n^*(f_1),\ldots, Z_n^*(f_k))$.
\end{proofof}

Next we present a useful technical lemma. It makes it possible
to replace  some of the assumptions of Theorem \ref{fidisconvtheo}
by sufficient conditions which are more restrictive  but often
simpler to verify.

\begin{lemma}\label{simplerfidicond}
\begin{enumerate}\item
If $Var(\Delta_n(f))=o(r_n v_n)$, then (C1) holds.
\item If $n v_n\to \infty$ and $\|f\|_\infty := \sup_{x\in
E_\cup}
  |f(x)|<\infty$, then (C1) and (C2) hold.
  \item If $r_nv_n \to 0$ and
  \begin{equation} \label{simplercovcond}
\frac{1}{r_nv_n}E\big( f(Y_n)g(Y_n)\big) \to c(f,g) \quad \forall\,
     f, g \in\FF,
\end{equation}
then (C3) holds.
\item If
\begin{equation} \label{simplerlindeberg}
E\Big( f(Y_n)^2\Ind{|f(Y_n)|>\eps\sqrt{nv_n}}\Big) = o(r_n v_n),
\quad \forall\,
    \eps>0, f\in\FF,
\end{equation}
then (C2) holds.
\item
If  $nv_n \to \infty$ and $\big(f(Y_n)^2\big)_{n\in\N}$ is
  uniformly integrable  under $P(\cdot)/(r_nv_n)$ for all $f\in\FF$, then (C2) holds.
\item If  $E\big(f(Y_n)^{2+\delta}\big)=O(r_nv_n)$
  for some $\delta>0$ and  all $f\in\FF$, then  $\big(f(Y_n)^2\big)_{n=1}^\infty$
  is uniformly integrable under $P(\cdot)/(r_nv_n)$ for all
$f\in\FF$.
\item If ($\widetilde{B3}$) holds, then $\lim_{k \to \infty} \limsup_{n \to \infty}
\frac{1}{r_nv_n} P\{L(Y_n) > k\}  = 0$ and the cluster lengths
$(L(Y_n))_{n\in\N}$ are tight under $P(\cdot\,|\, Y_n \neq 0)$.
\end{enumerate}
\end{lemma}
\begin{proof}

(i) The first equation in (C1) follows readily, the last one
by Chebyshev's inequality and the second one similarly using the inequality
 $$E\Big(|\Delta_n(f)-E\Delta_n(f)|\Ind{|\Delta_n(f)-E\Delta_n(f)|>\sqrt{nv_n}}\Big) \le \frac{Var(\Delta_n(f))}{\sqrt{nv_n}}.
 $$

 (ii) Under these conditions, (C2) obviously holds. Moreover, (C1)
  follows by (i), since
   $|\Delta_n(f)|\le
  2\|f\|_\infty \Ind{\Delta_n(f)\ne 0}$ implies
  \begin{eqnarray*}
    Var(\Delta_n(f)) & \le &
       E\Delta_n^2(f)\\
      & \le & 4 \|f\|_\infty^2 P\{\Delta_n(f)\ne 0\}\\
      & = & O\big( P\{ X_{n,i}\ne 0 \text{ \rm for some }
      r_n-l_n+1\le i\le r_n\}\big)\\
      & = & O(l_nv_n) \\
      & = & o(r_n v_n).
\end{eqnarray*}

(iii) By \eqref{simplercovcond}, $P\{Y_n \neq 0\} \leq
r_nv_n \to 0$ and the Cauchy-Schwarz inequality we have that
\begin{eqnarray} \label{meanbound}
\frac{1}{\sqrt{r_nv_n}}E|f(Y_n)| & = & \frac{1}{\sqrt{r_nv_n}}
E\big(|f(Y_n)|\Ind{Y_n \neq 0}\big)\nonumber \\
& \leq & \Big(\frac{1}{r_nv_n}E\big(f(Y_n)^2\big) P\{Y_n \neq 0\}
\Big)^{1/2} \to 0
  \end{eqnarray}
for $f \in \FF$. (C3) then follows readily from
\eqref{simplercovcond}.

(iv) By \eqref{simplerlindeberg}, for any $\epsilon > 0$,
\begin{eqnarray*}
E\Big(\Big(\frac{|f(Y_n)|}{\sqrt{nv_n}}\Big)^2\Big) & \leq & \epsilon^2
+ \frac{1}{nv_n}  E\big(f(Y_n)^2\boldsymbol{1}_{\{|f(Y_n)| >
\epsilon \sqrt{nv_n}\}}\big)\\
& = &\epsilon^2 +
o\Big(\frac{r_nv_n}{nv_n}\Big) = \epsilon^2 +o(1).
\end{eqnarray*}
Hence $Ef(Y_n) = o(\sqrt{nv_n})$, and (C2) then follows from
\eqref{simplerlindeberg}  by standard reasoning.

(v) By uniform integrability, $n/r_n \to \infty$ and Chebychev's
inequality,
$$
P\{|f(Y_n)| > \epsilon \sqrt{nv_n} \}  \leq \frac{ E\big(f(Y_n)^2
 \big)/(r_nv_n)}{ \epsilon^2 n/r_n} \to 0.
$$
Using uniform integrability again, it follows that
$E\big(f(Y_n)^2\boldsymbol{1}_{\{|f(Y_n)| > \epsilon \sqrt{nv_n}\}}
\big)$ $(r_nv_n)^{-1} \to 0$, so that \eqref{simplerlindeberg} is
satisfied. The result then follows from part (iv).

(vi)  is a well known fact.

(vii)
 Since, by stationarity,
\begin{eqnarray*}
\frac 1{r_nv_n} P\{L(Y_n)>k\} & \le & \frac 1{r_nv_n} \sum_{i=1}^{r_n-k} \sum_{j=i+k}^{r_n} P(X_{n,j}\ne 0|X_{n,i}\ne 0) P\{X_{n,i}\ne 0\} \\
& \le & \sum_{j=k}^{r_n} s_n(j),
\end{eqnarray*}
the assertion  follows readily from ($\widetilde{\text{B3}}$).
\end{proof}

\begin{proofof} {\sc Corollary  \ref{covcondcorol}.}\rm\quad
The first assertion follows if we prove that (C3') implies (C3).
However, using that $|E\big(f(Y_n)g(Y_n)\boldsymbol{1}_{\{L(Y_n) > k
\}}\big)| \leq$\\
 $\Big(E\big(f(Y_n)^2\boldsymbol{1}_{\{L(Y_n) > k \}}\big)
  E\big(g(Y_n)^2\boldsymbol{1}_{\{L(Y_n) > k \}}\big)\Big)^{1/2}$ it follows from
     \eqref{eq:C3.1'} and \eqref{eq:C3.2'}  that
\begin{eqnarray*}
  \lefteqn{\frac{1}{r_nv_n}E\big(f(Y_n)g(Y_n)\big)}\\
   &=&
 \frac{1}{r_nv_n}E\big(f(Y_n)g(Y_n)\boldsymbol{1}_{\{L(Y_n) \leq k \}}\big) +
 \frac{1}{r_nv_n}E\big(f(Y_n)g(Y_n)\boldsymbol{1}_{\{L(Y_n) > k \}}\big)\\
 &=& c_k(f,g) + R'_{n,k},
\end{eqnarray*}
with $\lim_{k \to \infty} \limsup_{n \to \infty} R'_{n,k} =0$.  A
standard subsequence argument then shows that $c(f,g) := \lim_{k \to
\infty} c_k(f,g)$ exits, and that
$$
\lim_{n\to\infty} \frac{1}{r_nv_n} E\big(f(Y_n)g(Y_n)\big) = c(f,g).
$$
By Lemma \ref{simplerfidicond} (iii) it then follows that  (C3) holds.

Now suppose instead  that (B1), (B2), (B3), (C1), and (C3'') hold.
The assumption (C2) then follows from Lemma \ref{simplerfidicond}
(v),  and hence only (C3) remains to be established.
 By Lemma \ref{lemma:clusterconv} (ii) and (iii),
$\theta_n=P\{Y_n\ne 0\}/(r_nv_n)\to \theta>0$ and $P^{(fg)(Y_n)\mid
Y_n\ne 0}$ converges weakly to $\mu_{fg,W}$. Thus, the uniform
integrability of  $(fg)(Y_n)$ under $P(\cdot)/(r_nv_n)$ is
equivalent to the uniform integrability under $P(Y_n \neq 0)$ so
that
\begin{eqnarray*}
 \frac 1{r_n v_n} E\big(f(Y_n)g(Y_n)\big)  &=&
\frac{P(Y_n \neq 0)}{r_nv_n} E\big(f(Y_n)g(Y_n) \mid Y_n \neq 0 \big)\\
 &\to& \theta\int x \, \mu_{fg,W}(dx) = E\big((f g)\big(W
\big)-(fg)\big(W^{(2; \infty)}\big)\big).
\end{eqnarray*}
It then follows from Lemma \ref{simplerfidicond} (iii) that (C3) holds with
$c(f,g)$ given by \eqref{thetacovconv}.
\end{proofof}

\begin{proofof} {\sc Lemma \ref{lemma:clusterconv}.}\rm\quad
 Again let $M_{n,s}^t := \sum_{i=s+1}^t \boldsymbol{1}_{\{X_{n,i}\ne 0\}}$ denote the
  number of non-vanishing observations in the time interval
  from $s+1$ to $t$. Then
  \begin{equation} \label{S6cond}
   \limsup_{n\to\infty} P\big( M_{n,l}^{r_n}\ne 0\mid X_{n,1}\ne
  0\big) \le \limsup_{n\to\infty} \big( \beta_{n,l} +
  r_nv_n\big)  \to 0
  \end{equation}
   as $l\to\infty$, by (B3) and $r_nv_n \to 0$. Hence, the
   analog to Condition (2) of Segers (2003) holds and one may
   conclude the assertions (i) and (ii) by essentially the same arguments as
   given for the proofs of Theorem 1 (with $t_n=r_n$), Corollary 2 and  Theorem 3 (i) there.

   The proof of (iii) also follows the ideas used in the proof of
   Theorem 3 (ii) in that paper. Nevertheless, we give more
   details, since we want to avoid working with the space
   $\mathbb{A}$ of sequences with almost all terms equal to 0, that
   was introduced by Segers (2003). Moreover, in this proof we replace
   assumption \eqref{fcontcond} in condition
   (C3.1'') by the weaker assumptions \eqref{fcontcond1} and
   \eqref{fcontcond2}.

   We first consider a bounded cluster functional $g$ such that $D_{g,m,I}\subset D_{f,m,I}$ for all
   $m\in\N$ and $I\subset\{1,\ldots,m\}$. The result for $f$ itself will then follow easily.
   Let $k\in\N$ be arbitrary and as before let $\|\cdot -
   \cdot\|_{TV}$ denote the total variation distance between two measures.
   By \eqref{S6cond}, for all $\eps>0$ there exists $l>k$ such that
   for sufficiently large $n$ and $X_n^{(k)}=(X_{n,i})_{1\le i\le
   k}$
   \begin{eqnarray}
     \lefteqn{\Big\| P\big( X_n^{(k)}\in\cdot,\;
     M_{n,k}^{r_n}=0\mid X_{n,1}\ne 0 \big) - P\big( X_n^{(k)}\in\cdot,\;
     M_{n,k}^{l}=0\mid X_{n,1}\ne 0 \big)\Big\|_{TV} } \nonumber\\
       & \le & P\big( M_{n,l}^{r_n}\ne 0\mid X_{n,1}\ne 0\big)\nonumber \hspace*{6cm}\\
       & \le & \eps \label{tvdist}
   \end{eqnarray}
   and, by \eqref{eq:finitew},
   \begin{eqnarray}
     \lefteqn{\Big\| P\big\{ W^{(k)}\in\cdot,\;
     W^{(k+1;\infty)}=0\big\} - P\big\{ W^{(k)}\in\cdot,\;
     W^{(k+1;l)}=0\big\} \Big\|_{TV} } \nonumber\\
       & \le & P\{ W_i\ne 0 \text{ for some } i>l\} \nonumber \hspace*{4cm}\\
   \label{tvdist2}    & \le & \eps.
   \end{eqnarray}
  Recall the definition of the sets $N_{k,I}$ for $I\subset \{1,\ldots,k\}$ from Remark
  \ref{condC3ppgen}.
  Since, according to assumption (C3.1''), the substochastic measures
  $P\big( X_n^{(k)}\in \cdot,\; X_n^{(k)}\in N_{k,I},\; M_{n,k}^l=0\mid
  X_{n,1}\ne 0\big)$ converge weakly to the substochastic measure
  $P\big\{ W^{(k)}\in \cdot,\; W^{(k)}\in N_{k,I},\;W^{(k+1;l)}=0\big\}$,
  it follows form \eqref{tvdist} and \eqref{tvdist2} that, for all $k\in\N$, and all subsets
  $I\subset \{1,\ldots,k\}$,
  \begin{eqnarray}
    \lefteqn{P\big( X_n^{(k)}\in \cdot,\; X_n^{(k)}\in N_{k,I},\; M_{n,k}^{r_n}=0\mid
    X_{n,1}\ne 0\big)}\nonumber \\
    \label{substochweakconv} & \to & P\big\{ W^{(k)}\in \cdot,\; W^{(k)}\in N_{k,I},\; W^{(k+1;\infty)}=0\big\}
  \end{eqnarray}
  weakly.

  By assertion (i) we have
   \begin{equation}   \label{theo1analog}
     E\big( g(Y_{n}) \mid Y_{n}\ne 0\big) = \frac 1{\theta_n} E\Big(
     g( X_n^{(r_n)}) - g( X_n^{(2;r_n)})\mid X_{n,1}\ne 0\Big) + o(1).
   \end{equation}
   Again by \eqref{tvdist} and the definition of a cluster
   functional,
   \begin{eqnarray}
     \lefteqn{ \Big| E\Big(
     g( X_n^{(r_n)}) - g( X_n^{(2;r_n)})\mid X_{n,1}\ne 0\Big)} \nonumber \\
      & & - E\Big(
     \big(g( X_n^{(l)}) - g( X_n^{(2;l)})\big) \boldsymbol{1}_{\{M_{n,l}^{r_n}=0\}}
     \mid X_{n,1}\ne 0\Big)\Big|  \; \le
     \; 2\eps \|g\|_\infty.    \label{gfinapprox}
   \end{eqnarray}
   In view of \eqref{substochweakconv} (with $k=l$), for all $I\subset\{1,\ldots,l\}$, the continuous
   mapping theorem yields
   \begin{eqnarray*}
    \lefteqn{E\Big(
     g( X_n^{(l)}) \boldsymbol{1}_{\{X_n^{(l)}\in N_{l,I}\}}
     \boldsymbol{1}_{\{M_{n,l}^{r_n}=0\}}\mid X_{n,1}\ne 0\Big)}\\
     & \to & E \Big( g(W^{(l)}) \boldsymbol{1}_{\{W^{(l)}\in N_{l,I}\}}
      \boldsymbol{1}_{\{W^{(l+1;\infty)}=0\}}\Big),
   \end{eqnarray*}
   because the function $g|_{N_{l,I}}$ is bounded and continuous on the
   complement of the set $D_{f,l,I}$, which
   by \eqref{fcontcond1} is a null set under the
   limit measure in \eqref{substochweakconv}. Sum up these
   equation for all $I\subset\{1,\ldots,l\}$ and combine this with
   an analogous result for $g(X_n^{(2;l)})$ to obtain
   \begin{eqnarray}
     \lefteqn{E\Big( \big(g( X_n^{(l)})-g(X_n^{(2;l)})\big)
     \boldsymbol{1}_{\{M_{n,l}^{r_n}=0\}}\mid X_{n,1}\ne
     0\Big)}\nonumber\\
     \label{gfinconv} &\to & E \Big( \big( g(W^{(l)})-g(W^{(2;l)})\big)
      \boldsymbol{1}_{\{W^{(l+1;\infty)}=0\}}\Big).
   \end{eqnarray}

   Combining \eqref{tvdist2}, \eqref{theo1analog}--\eqref{gfinconv}
   and $\theta_n\to\theta>0$, one arrives at
   \begin{equation}   \label{mainconvstep}
     E\big(g(Y_{n})\mid Y_{n}\ne 0\big) \to \frac 1\theta
     E\big(g(W)-g(W^{(2;\infty)})\big).
    \end{equation}
    Now, if $f$ is an arbitrary cluster functional satisfying the
   conditions of the proposition and $h:\R\to\R$ is continuous and
   bounded, then an application of \eqref{mainconvstep} with
   $g=h\circ f$ yields assertion (iii).
\end{proofof}

\begin{proofof} {\sc Corollary \ref{bddfunctfidicond}.}\rm\quad
  This is immediate from Corollary  \ref{covcondcorol} and Lemma
  \ref{simplerfidicond} (ii).
\end{proofof}

\begin{proofof} {\sc Theorem \ref{brackettheo}.}\rm\quad
  The processes $Z_n$ are asymptotically tight if the analogous sums
  over the even numbered and over the odd numbered blocks
  \begin{eqnarray} 
   \frac 1{\sqrt{n v_n}} & & \sum_{j=1}^{\floor{m_n/2}} \big( f(Y_{n,2j})-
E f(Y_{n,2j})\big) \quad \text{and} \nonumber \\
 \frac 1{\sqrt{n v_n}} & &
\sum_{j=1}^{\ceil{m_n/2}} \big( f(Y_{n,2j-1})- E f(Y_{n,2j-1})\big)  \label{evenoddblocks}
   \end{eqnarray}
 are asymptotically tight. In view of \eqref{totvardist}, the first
 expression is asymptotically tight if and only if the analogous
 expression with independent blocks, i.e.
  \begin{equation}  \label{totalvardist2}
   \frac 1{\sqrt{n v_n}} \sum_{j=1}^{\floor{m_n/2}} \big( f(Y_{n,2j}^*)-
E f(Y_{n,2j}^*)\big)
  \end{equation}
  is asymptotically tight, which follows from Theorem 2.11.9 of van
  der Vaart and Wellner (1996) applied with $Z_{ni}(f)=f(Y_{n,2i})$
  (and $m_n$ replaced with $\floor{m_n/2}$). Observe that for a sequence of
  monotonically increasing positive functions $T_n(\delta)$ the
  convergence of $T_n(\delta_n)$ to 0 for all sequences
  $\delta_n\downarrow 0$ is equivalent to $\lim_{\delta\downarrow 0}
  \limsup_{n\to\infty} T_n(\delta)=0$, so that the last two
  displayed conditions in Theorem 2.11.9 of van
  der Vaart and Wellner (1996) can be reformulated as (D3) and (D4), respectively.
  The proof of  tightness of the
  sum over the blocks with odd numbers is the same.
\end{proofof}

\begin{proofof} {\sc Remark \ref{condrem1} \rm(ii).}\rm\quad
 By the Cauchy-Schwarz inequality
    \begin{eqnarray*}
      \lefteqn{E^* \Big( F(Y_{n}) \Ind{F(Y_n)>\eps\sqrt{nv_n}}\Big)}\\
       &
      \le & \bigg( E^* \Big( F^2(Y_n)
      \Ind{F(Y_n)>\eps\sqrt{nv_n}}\Big) \cdot E^*
      \Ind{F(Y_n)>\eps\sqrt{nv_n}}\bigg)^{1/2}\\
      & \le & \bigg( \frac{\big(E^* \big( F^2(Y_n)
      \Ind{F(Y_n)>\eps\sqrt{nv_n}}\big)\big)^2}{\eps^2 n
      v_n}\bigg)^{1/2}\\
      & = & o\bigg( \frac{(r_nv_n)^2}{nv_n}\bigg)^{1/2}\\
      & = & o\big(r_n \sqrt{v_n/n}\big),
    \end{eqnarray*}
so (D2) holds. Further, (D2') implies \eqref{simplerlindeberg}, and
hence (C2) follows from Lemma \ref{simplerfidicond} (iv).

Next, suppose $E^*F^{2+\delta}(Y_n) = O(r_nv_n)$ and $nv_n \to
\infty$. Then
    \begin{eqnarray*}
      \lefteqn{E^*\Big( F^2(Y_n)
      \Ind{F(Y_n)>\eps\sqrt{nv_n}} \Big)}\\
       & \le & \big( E^*
       F^{2+\delta}(Y_n)\big)^{2/(2+\delta)}\cdot
       \big( E^*\Ind{F(Y_n)>\eps\sqrt{nv_n}}
       \big)^{1-2/(2+\delta)} \\
       & = & O\big((r_nv_n)^{2/(2+\delta)}\big) \cdot \bigg(
       \frac{E^*
       F^{2+\delta}(Y_n)}{\big(\eps\sqrt{nv_n}\big)^{2+\delta}}\bigg)^{1-2/(2+\delta)} \\
       & = & O\big( r_n v_n (n v_n)^{-\delta}\big)\\
       & = & o(r_n v_n),
    \end{eqnarray*}
    so that (D2') holds.
\end{proofof}

\begin{proofof} {\sc Theorem \ref{condrem2}.}\rm\quad
   First assume (D6) holds. Using the triangle inequality, it is easily seen that $Z_n$ is
   asymptotically equicontinuous if both terms given in
   \eqref{evenoddblocks} are asymptotically equicontinuous.
   Further, by \eqref{totvardist}, the
   first term is asymptotically equicontinuous if and only if
   \eqref{totalvardist2} is asymptotically equicontinuous. However,
   asymptotic equicontinuity of \eqref{totalvardist2}
   follows from Theorem 2.11.1 of van der Vaart and
   Wellner (1996). To see this, note that
   (D6) implies the analogous random entropy
   condition for the sums over the even numbered blocks, because
   the corresponding random semi-metric is smaller for these sums.

   If $m_n$ is even, then the second term in \eqref{evenoddblocks}
   has the same distribution as the first one, while for $m_n$ odd
   with probability greater than or equal to $1-r_nv_n\to 1$
   the additional summand $(nv_n)^{-1/2} (f(Y_{n,m_n})- E
   f(Y_{n,m_n}))$ equals $-(nv_n)^{-1/2}E   f(Y_{n,m_n})$, which
   tends to 0 uniformly for $f\in\FF$ (cf.\ \eqref{meanbound}).
   This proves the first assertion of the theorem.
   Theorem \ref{fidisconvtheo} then yields the convergence of $Z_n$,
   because the Lindeberg condition (C2) follows from (D2) (see
   Remark \ref{condrem1} (ii)).

   Next, to see that (D6') implies (D6), check that the random semi-metric $d_n$ can be represented as
  $d_n=(m_n/(nv_n))^{1/2}\cdot d_Q $ with the (random) probability measure
  $Q=m_n^{-1} \sum_{j=1}^{m_n} \eps_{Y_{n,j}^*}$, and hence
  $N(\eps,\FF,d_n)=N(\eps (nv_n/m_n)^{1/2},\FF,$ $d_Q)$. If $\int F^2\, dQ=0$, then
  $d_n(f,g)=0$ for all $f,g\in\FF$ and the integral in (D6') vanishes. Otherwise, for all
  $\eta>0$ there exists a $\tau>0$ such that for sufficiently large
  $n$
  $$ P\big\{{\textstyle (\int F^2\, dQ)^{1/2}}>\tau (nv_n/m_n)^{1/2}\big\} \le \frac{E
  F^2(Y_{n,1})}{\tau^2
  n v_n/m_n} \le \eta,
  $$
since $EF^2(Y_n) = O(r_nv_n)$, and thus with probability larger than $1-\eta$
  \begin{eqnarray*}
    \int_0^\delta \sqrt{\log N(\eps,\FF,d_n)}\, d\eps
      & = & \tau \int_0^{\delta/\tau} \sqrt{\log N(\eps\tau,\FF,d_n)}\,
      d\eps\\
      & \le & \tau \int_0^{\delta/\tau} \sup_{Q\in\mathcal{Q}}
       \sqrt{\log N(\eps{\textstyle (\int F^2dQ)^{1/2}},\FF,d_Q)}\,
      d\eps\\
      & \to & 0
  \end{eqnarray*}
  as $\delta\downarrow 0$, under (D6').
\end{proofof}

\begin{proofof} {\sc Corollary \ref{gentailcorol}.}\rm\quad
  Condition (D1) is satisfied since $F(x_1,\ldots,x_k)$ $\le \sum_{i=1}^k
  \phi_{\max}(x_i)$ and since $\phi_{max}$ is assumed to be
  measurable and bounded.
  Similarly, Condition (D2') follows from
  $F(Y_{n})\le r_n\|\phi_{\max}\|_\infty $, since $r_n = o(\sqrt{nv_n})$ by assumption.

  By Lemma \ref{simplerfidicond} (i), Assumption (C1) follows if we show that $Var(\Delta_n(f))=o(r_n
  v_n)$. Now,
  \begin{eqnarray*}
    E\Big( \sum_{i=1}^{r_n} \Ind{X_{n,i}\ne 0} \Big)^2 & \ge & E
    \sum_{j=1}^{\floor{r_n/l_n}} \Big( \sum_{i=1}^{l_n}
    \Ind{X_{n,(j-1)l_n+i}\ne 0} \Big)^2\\
    & = & \floor{r_n/l_n} E\Big( \sum_{i=1}^{l_n} \Ind{X_{n,i}\ne 0} \Big)^2
  \end{eqnarray*}
  by the row-wise stationarity, and consequently by
  \eqref{clustersizebound} and $l_n=o(r_n)$
  \begin{eqnarray*}
    E \big(\Delta_{n,1}^2(f)\big) & \le &  E\Big( \sum_{i=1}^{l_n} \phi_{\max}(X_{n,i})
    \Big)^2\\
    & \le & \|\phi_{\max}\|_\infty^2  E\Big( \sum_{i=1}^{l_n} \Ind{X_{n,i}\ne 0}
    \Big)^2 \\
    & = & O\Big(\frac{l_n}{r_n} r_n v_n\Big)\\
    & = & o(r_nv_n).
  \end{eqnarray*}

  Further, \eqref{gentailarrayapprox} follows from
  \begin{eqnarray*}
    \lefteqn{E^* \Big( \sup_{\phi\in\Phi} \frac 1{\sqrt{nv_n}} \Big|
    \sum_{i=r_nm_n+1}^n \big(
    \phi(X_{n,i})-E\phi(X_{n,i})\big)\Big|\Big)^2}\\
    & \le & E \Big( \frac 2{\sqrt{nv_n}} \|\phi_{\max}\|_\infty
    \sum_{i=r_nm_n+1}^n \Ind{X_{n,i}\ne 0}\Big)^2\\
    & = & \frac{4  \|\phi_{\max}\|^2_\infty}{nv_n} \cdot r_n v_n\\
    & \to & 0
  \end{eqnarray*}

    Therefore, the remaining assertions follow from Theorems \ref{brackettheo} and
    \ref{condrem2} and Remark \ref{condrem1} (i) and (ii).
\end{proofof}

\begin{proofof} {\sc Remark \ref{covrem} \rm (i).}\rm\quad
Since
\begin{eqnarray*}
\lefteqn{\frac{1}{r_nv_n} E \big(g_\phi(Y_n)^2
\boldsymbol{1}_{\{L(Y_n)
> k \}} \big)}\\
 &\leq& ||{\phi}||_\infty \frac{1}{r_nv_n} E \Big(
\Big( \sum_{i=1}^{r_n}\Ind{X_{n,i}\ne
0}\Big)^2\boldsymbol{1}_{\{L(Y_n)
> k\}} \Big)\\
&\leq&
 ||{\phi}||_\infty \bigg(\frac{1}{r_nv_n}E \Big(\Big(\sum_{i=1}^{r_n}\Ind{X_{n,i}\ne
0}\Big)^{2+\delta}\Big)\bigg)^{\frac{2}{2+\delta}}
\Big(\frac{1}{r_nv_n}P\{L(Y_n) > k\}\Big)^{\frac{\delta}{2+\delta}}
\end{eqnarray*}
 the first part \eqref{eq:C3.1'} of (C3') follows from \eqref{eq:ltightness}
 and \eqref{eq:twoplusdeltabound}, since $\phi$ is assumed to be bounded.
Next,
\begin{eqnarray} \label{eq:gphicovconv}
  \lefteqn{\frac{1}{r_nv_n}E\big(g_\phi(Y_n) g_\psi(Y_n)   \boldsymbol{1}_{\{L(Y_n) \leq k
  \}}\big)} \nonumber\\
   &=&
\frac{1}{r_nv_n}\sum_{i,j \in\{1,\ldots,r_n\}, |i-j| \leq k -1}
E\big( \phi(X_{n,i}) \psi(X_{n,j})\boldsymbol{1}_{\{L(Y_n) \leq k\}} \big) \nonumber \\
&=& \frac{1}{v_n}E\big( \phi(X_{n,1}) \psi(X_{n,1})\big) \\
 &&{ } + \sum_{i =1}^{k-1} \frac{r_n - i}{r_n} \frac{1}{v_n}
\big( E(\phi(X_{n,1}) \psi(X_{n,i+1})) + E(\psi(X_{n,1})
\phi(X_{n,i+1})\big)+ R_{n,k},\hspace*{-2cm} \nonumber
\end{eqnarray}
with
\begin{eqnarray*}
|R_{n,k}| &=& \frac{1}{r_nv_n}\Big|\sum_{i,j \in\{1,\ldots,r_n\},
|i-j| \leq k -1}
E\big( \phi(X_{n,i}) \psi(X_{n,j})\boldsymbol{1}_{\{L(Y_n) > k\}} \big)\Big| \\
&\leq& ||{\phi}||_\infty ||\psi||_\infty \frac{1}{r_nv_n} E\Big(
\Big(\sum_{i=1}^{r_n}\Ind{X_{n,i}\ne
0}\Big)^2\boldsymbol{1}_{\{L(Y_n)
> k\}} \Big).
\end{eqnarray*}
It then follows as above that $\lim_{k \to \infty} \limsup_{n \to
\infty}|R_{n,k}| = 0$, and hence the assumption \eqref{eq:C3.2'} of
(C3') can be seen to be satisfied, with $c$ given by
\eqref{eq:dsum}.

Furthermore, if (B1), ($\widetilde{\text{B3}}$) and \eqref{eq:phicovconv0} are fulfilled, then by stationarity
\begin{align*}
  \frac 1{r_nv_n}& \Cov\big(g_\phi(Y_n),g_\psi(Y_n)\big) \\
  & = \frac 1{r_nv_n} E\big(g_\phi(Y_n)g_\psi(Y_n)\big) + O(r_nv_n)\\
  & = \frac 1{v_n} E(\phi(X_{n,1})\psi(X_{n,1})) +o(1) \\
  & \hspace*{0.5cm} + \sum_{k=1}^{r_n-1} \frac{1-k/r_n}{v_n} \big(E(\phi(X_{n,1})\psi(X_{n,k+1})) +E(\phi(X_{n,k+1})\psi(X_{n,1}))\big)\\
  & \to c(g_\phi,g_\psi).
\end{align*}
In the last step we may apply Pratt's lemma (Pratt, 1960), because $\phi$ and $\psi$ are bounded and ($\widetilde{\text{B3}}$) holds.
\end{proofof}

\begin{proofof} {\sc Corollary \ref{gentailcorol2}.}\rm\quad
  Clearly \eqref{clustersumbound} implies
  \eqref{ljapcond} and hence also (D2'). Moreover, \eqref{clustersumbound} implies that
  \begin{eqnarray*}
    E\Big( \sum_{i=1}^{r_n} \phi_{\max}(X_{n,i})\Big)^2 & \le &
      E\Big( \sum_{i=1}^{r_n} \phi_{\max}(X_{n,i})\Big)^{2+\delta}\\
      & & \hspace*{1cm} +
      P\Big\{ 0<\sum_{i=1}^{r_n} \phi_{\max}(X_{n,i}) \le 1\Big\}\\
      & = & O(r_n v_n).
  \end{eqnarray*}
  Hence, similar arguments as used in the proof of Corollary
  \ref{gentailcorol} show that $(Z_n(g_\phi))_{\phi\in\Phi}$
  converges weakly to a Gaussian process. Finally,
  \eqref{gentailarrayapprox} and thus the convergence of $(\tilde
  Z_n(\phi))_{\phi\in\Phi}$ follows from
  \begin{eqnarray*}
    \lefteqn{E^* \Big( \sup_{\phi\in\Phi} \frac 1{\sqrt{nv_n}} \Big|
    \sum_{i=r_nm_n+1}^n \big(
    \phi(X_{n,i})-E\phi(X_{n,i})\big)\Big|\Big)^2} \\
    & \le & E \Big( \frac 1{\sqrt{nv_n}}
    \sum_{i=1}^{r_n} \big(\phi_{\max}(X_{n,i})+E \phi_{\max}(X_{n,i})\big)\Big)^2\\
    & \le & \frac{4 }{nv_n}  E \Big( \sum_{i=1}^{r_n} \phi_{\max}(X_{n,i})\Big)^2\\
    & = & O(r_n/n)\\
    & \to & 0.
  \end{eqnarray*}
\end{proofof}

\begin{proofof} {\sc Corollary \ref{cor:indicators}.}\rm \;
(i) The index set
  $ \CC:=\{ C_{t_1,\ldots,  t_k} \mid t_1,\ldots,t_k\in [0,1]\} $
   equipped with the metric
  $ \rho_\CC(\boldsymbol{1}_{C_{s_1,\ldots, s_k}},\boldsymbol{1}_{C_{t_1,\ldots,  t_k}})$ $ :=
  \max_{1\le l\le k} |s_l-t_l|
  $
  is totally bounded. The same holds for $ \DD:=\{ D_{t_1,\ldots,  t_k} \mid t_1,\ldots,t_k\in [0,1]\} $.

  In view of the discussion preceding Corollary \ref{cor:indicators},
  the assertions follow from Theorem \ref{condrem2} combined  with
  Corollary \ref{bddfunctfidicond} if we verify
  condition (D5) and that the index sets $\CC$ and $\DD$ are
  VC-classes. Condition (D5) is satisfied since all processes under
consideration are separable.

  That $\CC$ is a VC-class may be established by observing that
  $ C_{t_1,\ldots,  t_k} = \psi^{-1}\big( \times_{l=1}^k
  (t_l,\infty)\big)$ with
  \begin{align*}
   \psi:\, & \R_\cup \to \R^k\\
     & (x_1,\ldots,x_m)\mapsto \left\{
  \begin{array}{l@{\quad}l}
    (x_{j},\ldots,x_{j+k-1}) & \text{if } j=\min\{i\mid x_i \neq 0\}\le
    m-k+1\\
    (0,\ldots,0) & \text{else}.
  \end{array}
  \right.
 \end{align*}
 Since $\{ \times_{l=1}^k  (t_l,\infty) \mid t_1,\ldots, t_k\ge 0\}$
 is known to be a VC-class (cf.\ van der Vaart and Wellner, 1996,
 Example 2.6.1), $\CC$ is a VC-class, too (van der Vaart and Wellner,
 1996, Lemma 2.6.17 (v)).

 The sets $\DD_j := \{E_{j,t}\mid t\ge 0\}$ are linearly
  ordered (i.e., $E_{j,s}\subset E_{j,t}$ if $s>t$) and hence they are
  VC-classes, and hence so is
  $$ \DD = \DD_1\sqcap \DD_2\sqcap\cdots\sqcap \DD_k = \big\{
  \bigcap_{j=1}^k E_j \mid E_j\in\DD_j\big\}
  $$
  (van der Vaart and Wellner, 1996, Lemma 2.6.17 (ii)).
\smallskip

(ii) By the results of Segers (2003), Condition (C3.1'') is satisfied
in the weaker version discussed in Remark \ref{condC3ppgen}, because
the limit rv's are continuous on $(0,\infty)$ and the discontinuity
sets have Lebesgue measure 0. Hence the assertions follow by part
(i), if the asymptotic equicontinuity condition (D3) can be shown.

For this, first note that $C_{s_1,\ldots,
 s_k}\triangle$ $ C_{t_1,\ldots, t_k} \subset \big\{
 (x_1,\ldots,x_m)\in E_\cup\mid m\in\N,$ $ \exists\, \, 0\le j\le m-k, 1\le l\le k\; \forall 1\le i\le j:
 x_i=0, x_{j+l} \in (\min(s_l,t_l),\max(s_l,t_l)]\big\}$.
 Thus, Lemma \ref{lemma:clusterconv} (i) and (ii) yield that
 \begin{eqnarray*}
   \lefteqn{ \frac 1{r_nv_n} P\{Y_{n}\in C_{s_1,\ldots,
 s_k}\triangle C_{t_1,\ldots, t_k} \}} \\
   & \le & \frac 1{r_nv_n\theta_n} P\big(X_n^{(r_n)}\in C_{s_1,\ldots,
 s_k}\triangle C_{t_1,\ldots, t_k}\mid X_{n,1}\ne 0\big) \cdot
 P\{Y_{n}\ne 0\}\\
   & & \hspace*{8cm} + o\Big(\frac{P\{Y_{n}\ne 0\}}{r_nv_n}\Big)\\
  & = & P\big(X_n^{(r_n)}\in C_{s_1,\ldots,
 s_k}\triangle C_{t_1,\ldots, t_k}\mid X_{n,1}\ne 0\big)  +o(1)\\
 & \le & \sum_{l=1}^k P\big(X_{n,l}\in (\min(s_l,t_l),\max(s_l,t_l)]
 \mid X_{n,1}\ne 0\big) +o(1)\\
 & \le & \sum_{l=1}^k P\big(X_{n,l}\in (\min(s_l,t_l),\max(s_l,t_l)]
 \mid X_{n,l}\ne 0\big)\cdot \frac{P\{X_{n,l}\ne 0\}}{P\{X_{n,1}\ne 0\}} +o(1)\\
 & = &  \sum_{l=1}^k |t_l-s_l| +o(1),
 \end{eqnarray*}
 where the term $o(1)$ tends to 0 uniformly for all $s_1,\ldots,
 s_k, t_1,\ldots, t_k\in[0,1]$. Now, (D3) follows
 immediately from the definition of $\rho_\CC$.

 To verify condition (D3) for the indicator functions describing the
 largest  order statistics in a cluster, note that
\begin{eqnarray*}
   \lefteqn{\bigcap_{j=1}^k E_{j,s_j} \triangle \bigcap_{j=1}^k
   E_{j,t_j}}\\
    & \subset & \Big\{ (x_1,\ldots,x_m)\in E_\cup \mid m\in\N, \\
    & & \hspace*{0.1cm} 
    \sum_{i=1}^m
  \boldsymbol{1}_{(\min(s_j,t_j),1]}(x_i)\ge j, \sum_{i=1}^m \boldsymbol{1}_{(\max(s_j,t_j),1]}(x_i)<
  j \mbox{ for some } 1\le j\le k\Big\}\\
    & \subset & \big\{ (x_1,\ldots,x_m)\in E_\cup \mid m\in\N,\\
        & & \hspace*{1cm} x_i\in(\min(s_j,t_j),\max(s_j,t_j)] \mbox{ for some } 1\le j\le k, 1\le i\le
  m\big\}.
  \end{eqnarray*}
  This implies
  \begin{eqnarray*}
    \lefteqn{\frac 1{r_n v_n} P\Big\{Y_{n}\in \bigcap_{j=1}^k E_{j,s_j} \triangle \bigcap_{j=1}^k
    E_{j,t_j}\Big\}}\\
    & \le & \sum_{j=1}^k P \big( X_{n,1} \in (\min(s_j,t_j),\max(s_j,t_j)]
    \mid X_{n,1}\ne 0\big) \\
    & = & \sum_{j=1}^k |t_j-s_j|,
  \end{eqnarray*}
  from which (D3) follows.
\end{proofof}

\begin{proofof} {\sc the result in Example \ref{ex:allvalues}.}\rm\quad
The convergence of the fidis of $\tilde{Z}_n$ to those of a Gaussian
process with covariance function \eqref{eq:allvaluescov} follows
from Corollary \ref{bddfunctfidicond} by the same arguments as in
the proof of Corollary \ref{cor:indicators} (ii).

In view of the discussion before Corollary \ref{cor:indicators}, the
proof will be completed by showing that the conditions (D3), (D5)
and (D6) of the asymptotic equicontinuity Theorem \ref{condrem2}
also are satisfied. The measurability condition (D5) holds since,
for fixed  $k$, the processes
 $(1_{C_{k,t_1,\ldots,t_k}})_{(t_1,\ldots,,t_k)\in[0,1]^k}$ are
 separable and a supremum of countably many suprema of separable
 processes are measurable.

We will use \eqref{ascontdiff} to verify that (D3) is satisfied for
the semi-metric
\begin{eqnarray*}
 \lefteqn{ \rho(\boldsymbol{1}_{C_{j, s_1, \ldots, s_j}},\boldsymbol{1}_{C_{k, t_1, \ldots, t_k}})}\\
 & := & \left\{
      \begin{array}{l@{\quad}l}
        P\{L(W) \in \{j,k\}\} & \text{if } j \neq k,\\
         P\{L(W) =k, W_i \in (s_i \wedge t_i, s_i \vee t_i]
         \; \text{for some} \; 1 \leq i \leq k\}     & \text{if } j=k.
      \end{array} \right.
 \end{eqnarray*}
    Now, $\FF = \{\boldsymbol{1}_{C_{k, t_1, \ldots ,t_k}} \mid k \geq
1, t_1, t_2, \ldots \in [0,1]\}$  is totally bounded with respect to
$\rho$. To see this, for $\epsilon > 0$ given, choose $0=a_{i,0} <
a_{i,1} < \ldots <a_{i,m_i} =1$ such that $P\{W_i \in (a_{i,j-1},
a_{i,j}]\} \leq \epsilon/k_\epsilon$ for $1 \leq i \leq k_\epsilon$
and $1 \leq j \leq m_i$, with $k_\epsilon$ chosen large enough to
make $P\{L(W) \geq k_\epsilon\} < \epsilon/2$. Then
\begin{eqnarray*}
\{\boldsymbol{1}_{C_{k, t_1, \ldots ,t_k}} \mid k \geq k_\epsilon\},
\;
  \{\boldsymbol{1}_{C_{j, t_1, \ldots ,t_j}}
   \mid t_i \in [a_{i, \ell_{i}-1}, a_{i, \ell_i}], \forall 1\leq i \leq j\},
\end{eqnarray*}
for $1 \leq j \leq k_\epsilon, \; 1 \leq \ell_i \leq m_i$ is a
finite cover of $\FF$ with diameter at most $\epsilon$.

By Lemma \ref{lemma:clusterconv}
\begin{equation} \label{eq:rhoconv1}
P(L(Y_n)=k \mid Y_n \neq 0) \to \frac{1}{\theta}\Big(P\{L(W)=k\} -
P\{L(W^{(2; \infty)}) =k\} \Big),
\end{equation}
and by Sheffe's Lemma the convergence is uniform in $k\in\N$. (Note
that, for $k\le l$, the cluster functional  $\boldsymbol{1}_{\{k\}}
\circ L$ is constant on all sets $N_{l,I}$ defined in Remark
\ref{condC3ppgen}.) Similarly,
\begin{eqnarray}  \label{eq:rhoconv2}
\lefteqn{P\big(L(Y_n)=k, (Y_n^c)_1 \leq t_1, \ldots , (Y_n^c)_k \leq
t_k
  \mid Y_n \neq 0\big) } \nonumber\\
  &\to& \frac{1}{\theta}\Big(P\big\{L(W)=k, W_1 \leq t_1,
  \ldots, W_k \leq t_k\big\} \nonumber\\
 && { } - P\big\{L(W^{(2;\infty)})=k,\,
  \big((W^{(2;\infty)})^c\big)_i\le t_i,\, \forall 1\le i\le k \big\} \Big),
\end{eqnarray}
and the convergence is uniform in $t_1, \ldots ,t_k$ for each fixed
$k$, because the right-hand side defines a continuous function.

For $\epsilon > 0$  let $\delta = \epsilon/2$ and consider $j, t_1,
\ldots ,t_j, \; k, t_1, \ldots ,t_k$  such that
$\rho(\boldsymbol{1}_{C_{j, s_1, \ldots, s_j}},\boldsymbol{1}_{C_{k,
t_1, \ldots, t_k}})$ $ < \delta$. Then for $j \neq k$ and $n$ large
\begin{eqnarray*}
\frac{1}{r_nv_n}P\big\{Y_n \in C_{j, s_1, \ldots, s_j} \Delta C_{k,
t_1, \ldots, t_k} \big\}
&\leq& \frac{1}{r_nv_n}P\big\{L(Y_n) \in \{j,k\}\big\}\\
&=& \theta_n P\big(L(Y_n) \in \{j,k\} \mid Y_n \neq 0\big) \leq
\epsilon,
\end{eqnarray*}
by \eqref{eq:rhoconv1}, Lemma \ref{lemma:clusterconv} (ii) and the
definition of $\rho$.

If instead $j=k \leq k_\epsilon$ then using \eqref{eq:rhoconv2}, for large $n$,
\begin{eqnarray*}
\lefteqn{ \frac{1}{r_nv_n}P\big\{Y_n \in C_{j, s_1, \ldots, s_j}
\Delta C_{k, t_1, \ldots, t_k} \big\}
 }\\
&= &  \theta_n P\big(L(Y_n)=k, \;  (Y_n^c)_i \in (s_i \wedge t_i,
s_i \vee t_i]  \text{ for some} \ 1 \leq i \leq k \mid Y_n \neq 0\big)\\
&\leq& \theta_n \Big(
 \frac{1}{\theta}P\big\{L(W)=k, \;  W_i \in (s_i \wedge t_i, s_i \vee t_i]
  \text{ for some} \ 1 \leq i \leq k \big\} + \frac{\epsilon}{4}\Big)\\
  & \le & \epsilon,\\
\end{eqnarray*}
again by Lemma \ref{lemma:clusterconv}  and the definition of
$\rho$.

 Finally, if $j=k>k_\epsilon$, then for
large $n$
\begin{align*}
\frac{1}{r_nv_n} & P\big(Y_n \in C_{j, s_1, \ldots, s_j} \Delta C_{k,
t_1, \ldots, t_k} \big)\\
 & \leq P(L(Y_n) = k \mid Y_n \neq 0) \leq
2P(L(W) > k_\epsilon) < \epsilon.
\end{align*}
This concludes the proof of \eqref{ascontdiff}, and hence also the proof of (D3).

For the proof of (D6), let $\CC_k = \{C_{j, t_1, \ldots ,t_j} \mid
1\leq j \leq k, t_1,  \ldots ,t_j \in [0,1]\}$ and $\FF_k =
\{\boldsymbol{1}_C \mid C \in \CC_k\}$ so that $\FF =
\bigcup_{k=1}^\infty \FF_k$. Define $\psi_k$ as the function which
maps $x \in E_\cup$ to the vector $(1, \ldots , 1)$ in $\R^{2k}$ if
$L(x) > k$ or $L(x)=0$ and which maps $x$ to the vector
$$
(1, \ldots 1, 0, 1, \ldots 1, x_1^c, \ldots , x_j^c, 0, \ldots , 0)
\in \R^{2k}
$$
if $1 \leq L(x):= j \leq k$. Here the first row of ones has $j-1$
entries and the second  row has $k-j$ entries, and hence the vector
ends with $k-j$ zeros, so that the first $k$ components encode the
length of the cluster core. With this definition it follows that
$$
C_{j, t_1, \ldots ,t_j} = \psi_k^{-1}(\R^{j-1} \times (- \infty, 0]
\times \R^{k-j} \times \times_{i=1}^j (-\infty, t_i] \times
\R^{k-j}).
$$
The left orthants $\times_{i=1}^{2k} (-\infty, x_i]$ form a VC-class
with index bounded  by $2k+1$ (van de Vaart and Wellner (1996,
Example 2.6.1)) and hence also $\CC_k$ is a VC-class with index
bounded by  $2k+1$ (Dudley (1999, Theorem 4.2.3)).  By van de Vaart
and Wellner (1996, Theorem 2.6.7), for all sufficiently small
$\epsilon$ and all $k\in\N$,  $\FF_k$ satisfies the  metric entropy
bound
\begin{equation}  \label{eq.Nbound}
N\big(\epsilon({\textstyle\int} F^2 dQ)^{1/2}, \FF_k, d_Q\big) \leq
C(2k+1)(16\e)^{2k+1} \epsilon^{-(4k+1)} \leq \epsilon^{-(6k+2)} ,
\end{equation}
with $C$ denoting a universal constant that does not depend on $k$
or $\epsilon$.

Let $L_{n,1} > L_{n,2} > \ldots L_{n,{m_n}}$ be the order statistics
in descending order  of the independent cluster lengths
$\big(L(Y_{n, j}^*)\big)_{j=1}^{m_n}$. Since the empirical
$L_2$-semi-metric $d_n$ satisfies
\begin{equation*}
\sup_{i,j > k} d_n^2 \big(\boldsymbol{1}_{C_{i, t_1, \ldots ,t_i}},
\boldsymbol{1}_{C_{j, s_1, \ldots ,s_j}} \big) \leq \frac{1}{nv_n}
\sum_{j=1}^{m_n} \Ind{L(Y_{n, j}^*) >k}
\end{equation*}
it follows that the squared diameter of the set
$$
\{C_{j, t_1, \ldots ,t_j} \mid j> L_{n, \lfloor \epsilon^2 n v_n
\rfloor}, t_1, \ldots ,t_j \in [0, 1] \}
$$
w.r.t.\ $d_n$ is bounded by
$$
\frac{1}{nv_n} \sum_{j=1}^{m_n} \Ind{L(Y_{n, j}^*) > L_{n, \lfloor
\epsilon^2 n v_n \rfloor}}
 \leq \frac{\lfloor \epsilon^2 n v_n \rfloor}{nv_n} \leq \epsilon^2.
$$

Reasoning as in the last part of the proof of Theorem
\ref{condrem2}, this  together with \eqref{eq.Nbound} shows that
(D6) follows if we prove that
\begin{equation}  \label{eq:orderstatint}
\lim_{\delta \downarrow 0} \limsup_{n \to \infty} P\Big\{
\int_0^\delta \sqrt{\log \epsilon^{-(6L_{n, \lfloor \epsilon^2 n v_n
\rfloor}+2)}}\: d \epsilon >\tau\Big\} =0,
\end{equation}
for all $\tau>0$. By a change of variables and H\"older's inequality
\begin{eqnarray*}
\lefteqn{\int_0^\delta \sqrt{\log \epsilon^{-(6L_{n, \lfloor
\epsilon^2 n v_n \rfloor}+2)}}\, d \epsilon}\\
 &\le&
 \sum_{j=1}^{\lceil \delta n v_n
\rceil} \sqrt{8L_{n,j}}
\int_{(j/(nv_n))^{1/2}}^{((j+1)/(nv_n))^{1/2}}
\sqrt{|\log \epsilon|}\: d \epsilon \\
& \le &  \frac 2{nv_n}  \sum_{j=1}^{\lceil \delta n v_n \rceil}
\sqrt{L_{n,j}}\cdot nv_n\int_{j/(nv_n)}^{(j+1)/(nv_n)}
\sqrt {\log \eta^{-1/2}} \eta^{-1/2}\, d\eta\\
 &\leq& \Big(\frac{1}{nv_n} \sum_{j=1}^{\lceil \delta
n v_n \rceil}
 L_{n,j}^{1 + \zeta}\Big)^{1/(2+2\zeta)} \Big(\frac{1}{nv_n} \sum_{j=1}^{\lceil \delta n v_n
\rceil} \\
& & \hspace*{0.8cm}\Big(nv_n \int_{j/(nv_n)}^{(j+1)/(nv_n)} \sqrt {\log
\eta^{-1/2}} \eta^{-1/2} \: d \eta \Big)^{(2+2\zeta)/(1+2\zeta)}
\Big)^{(1+2\zeta)/(2+2\zeta)}.
\end{eqnarray*}
Now,
$$
E\Big(\frac{1}{nv_n} \sum_{j=1}^{\lceil \delta n v_n \rceil}
L_{n,j}^{1 + \zeta}\Big) \leq E\Big(\frac{1}{nv_n} \sum_{j=1}^{m_n}
L_n(Y_{n,j}^*)^{1 + \zeta}\Big) \le E\big(L(Y_n)^{1+\zeta}\mid
Y_n\ne 0\big),
$$
which is bounded by \eqref{eq:lengthhalfbound}. Furthermore,
applying Liapunov's  inequality to the individual summands,
\begin{eqnarray*}
 \lefteqn{\frac{1}{nv_n} \sum_{j=1}^{\lceil \delta n v_n \rceil} \Big(nv_n
\int_{j/(nv_n)}^{(j+1)/(nv_n)} \sqrt {\log \eta^{-1/2}} \eta^{-1/2}
\: d
\eta \Big)^{(2+2\zeta)/(1+2\zeta)}}\\
 &\leq& \frac{1}{nv_n} \sum_{j=1}^{\lceil \delta n v_n \rceil} nv_n
 \int_{j/(nv_n)}^{(j+1)/(nv_n)}
 \Big(\frac{|\log \eta|}\eta\Big)^{(1+\zeta)/(1+2\zeta)}\: d \eta\\
&\le & \int_0^{2\delta} \Big(\frac{|\log
\eta|}\eta\Big)^{(1+\zeta)/(1+2\zeta)}\: d \eta\to 0,
\end{eqnarray*}
as $\delta \to 0$. Hence we have verified \eqref{eq:orderstatint}. This
concludes the proof of (D6).
\end{proofof}

\bigskip

\noindent {\large\bf References}
   \smallskip

\parskip1.2ex plus0.2ex minus0.2ex

\rueck
Bortot, P. and Tawn, J.A. (1998). Models for the extremes of Markov chains. {\em Biometrika}
{\bf 85}, 851-867.


\rueck
  Drees, H. (2000). Weighted approximations of tail processes for $\beta$--mixing
  random variables. {\em Ann.\ Appl.\ Probab.} {\bf 10}, 1274--1301.

\rueck
  Drees, H. (2002). Tail empirical processes under mixing
  conditions. In H.G.\ Dehling, T.\ Mikosch and M.\ S{\o}rensen (eds.),
  {\em Empirical Process Techniques for Dependent Data}, 325-342. Boston: Birkh\"auser.

\rueck
 Drees, H. (2003). Extreme quantile estimation for dependent data with applications to
finance. {\em Bernoulli} {\bf 9}, 617--657.

\rueck
 Drees, H. (2009). Smoothed blocks estimator of the extremal index.
 Preprint, University of Hamburg.

\rueck
Dudley, R. (1999). {\em Uniform central limit theorems.} Cambridge: Cambridge University Press.

\rueck
  Eberlein, E. (1984). Weak convergence of partial sums of absolutely regular
sequences. {\em Statist.\ Probab.\ Lett.} {\bf 2}, 291--293.

\rueck
   Einmahl, J.  (1997). Poisson and Gaussian approximation of weighted local
empirical processes. {\em  Stochastic Process.\  Appl.} {\bf 70},
31--58.

\rueck
  Fabian, V. (1970). On Uniform Convergence of Measures. {\em
  Probab.\ Theory Related Fields} {\bf 15}, 139--143.

\rueck
  Gin\'e, E., and  Mason, D. M. (2008). Uniform in Bandwidth Estimation
of Integral Functionals of the Density Function. {\em Scand.\ J.\
Statist.} {\bf 35}, 739--761.

\rueck
 Gin\'e, E., Mason, D. M., and Zaitsev, A. Y. (2003). The L1-norm
density estimator process. {\em Ann.\ Probab.} {\bf 31}, 719--768.

\rueck
  de Haan, L., and Ferreira, A. (2006). {\em Extreme Value Theory.}
  Springer.

\rueck Hahn, M. G. (1977). Conditions for sample-continuity and the
central limit theorem. {\em Ann.\ Probab.} {\bf 5}, 351--360.

\rueck
  Leadbetter, M.R. (1995). On high level exceedance modeling and tail inference. {\em
J.\ Statist.\ Plann.\ Inference} {\bf 45}, 247--260.

\rueck
  Leadbetter, M.R., and Rootz\'{e}n, H. (1993). On central limit theory for families of strongly
  mixing additive random functions.
  In: Stochastic processes: a festschrift in honour of
Gopinath Kallianpur (S.\ Cambanis et al., eds.), 211--223.
Springer.

\rueck
Pratt, J.\ W.\ (1960). On interschanging limits and integrals. {\em Ann.\ Math.\ Statist.} {\bf 31}, 74--77.

\rueck
  Petrov, V.V. (1975). {\em Sums of Independent Random
  Variables.} Springer.

\rueck
  Rootz\'en, H. (1995). The tail empirical process for
  stationary sequences. Preprint, Chalmers University Gothenburg.

\rueck
  Rootz\'en, H. (2009): Weak convergence of the tail empirical function
  for dependent sequences.  {\em Stochastic Process.\   Appl.} {\bf 119},
  468--490.

\rueck
  Rootz\'en, H., Leadbetter, M.R., and de Haan, L. (1990). Tail
  and quantile estimators for strongly mixing stationary
  processes. Report, Department
  of Statistics, University of North Carolina.

\rueck
  Rootz\'en, H., Leadbetter, M.R., and de Haan, L. (1998). On the distribution of tail
  array sums for strongly mixing stationary sequences. {\em Ann.\ Appl.\ Probab.} {\bf
  8}, 868--885.

\rueck
  Segers, J. (2003). Functionals of clusters of extremes. {\em Adv.\ in Appl.\ Probab.} {\bf
  35}, 1028--1045.

\rueck
 Sisson, S., and Coles, S. (2003). Modelling Dependence
Uncertainty in the Extremes of Markov Chains. {\em Extremes} {\bf
6}, 283--300.

\rueck
 Stott, P. A., Stone, D. A., and Allen, M. R. (2004). Human
contribution to the European heatwave of 2003. {\em Nature} {\bf
432},  610--613.

\rueck
  van der Vaart, A.W., and Wellner, J.A. (1996). {\em Weak Convergence
  and Empirical Processes.} Springer.

\rueck
 Wellner, J.A., and Zhang, Y. (2000). Two estimators of the mean of a
counting process with panel count data. {\em Ann.\ Statist.} {\bf
28}, 779--814.

\rueck
 Wellner, J.A., and Zhang, Y. (2007). Two likelihood-based
semiparametric estimation methods for panel count data with
covariates {\em Ann.\ Statist.} {\bf 35}, 2106--2142.

\rueck
  Yun, S. (2000). The distribution of cluster functionals of extreme events in a
  $d$th-order Markov chain. {\em J.\ Appl.\ Probab.} {\bf 37}, 29--44.

\end{document}